\documentclass[11pt]{article}

\RequirePackage[utf8]{inputenc}
\RequirePackage[OT1]{fontenc}
\RequirePackage{amsthm,amsmath,amsfonts}
\RequirePackage[numbers]{natbib}
\RequirePackage[colorlinks,citecolor=blue,urlcolor=blue]{hyperref}
\RequirePackage{authblk}

\newtheorem{Th}{Theorem}[section]
 
\newtheorem{Cy}[Th]{Corollary} 
\newtheorem{Lmm}[Th]{Lemma} 
 
\newtheorem{Defn}[Th]{Definition}
 
 \title{Effect of Stochastic Perturbations for Front Propagation in Kolmogorov Petrovskii Piscunov Equations}  
 
\author{John M. Noble \thanks{email address: {\tt noble@mimuw.edu.pl}}}
\affil{Mathematical Statistics\\ Institute of Applied Mathematics and Mechanics,\\  Faculty of Mathematics, Informatics and Mechanics, \\  University of Warsaw}

\date{}
 
\begin{document}
 \maketitle
 \begin{abstract}
 This article considers equations of Kolmogorov Petrovskii Piscunov type in one space dimension, with stochastic perturbation:
 
 \[ \left\{ \begin{array}{l} \partial_t u = \left ( \frac{\kappa}{2} u_{xx} + u(1-u) \right ) dt + \epsilon u  \partial_t \zeta \\
             u_0(x) = {\bf 1}_{(-\infty, - \frac{1}{N}\log 2)}(x) + \frac{1}{2}e^{-Nx}{\bf 1}_{[-\frac{1}{N} \log 2, +\infty)}(x) 
            \end{array}\right. 
\]
 
\noindent where the stochastic differential is taken in the sense of Itô and $\zeta$ is a Gaussian random field satisfying $\mathbb{E} \left [ \zeta \right ] = 0$ and $\mathbb{E} \left [ \zeta(s,x)\zeta(t,y) \right ] = (s \wedge t) \Gamma (x-y)$. Two situations are considered: firstly, $\zeta$ is simply a standard Wiener process (i.e. $\Gamma \equiv 1$): secondly,  $\Gamma \in C^\infty (\mathbb{R})$ with $\int_{-\infty}^\infty |\Gamma(z)| dz < +\infty$.  

The results are as follows: in the first situation (standard Wiener process: i.e. $\Gamma(x) \equiv 1$), there is a non-degenerate travelling wave front if and only if $\frac{\epsilon^2}{2} < 1$, with asymptotic wave speed  $\max\left(\sqrt{2\kappa (1 - \frac{\epsilon^2}{2})}, \frac{1}{N}(1 - \frac{\epsilon^2}{2}) + \frac{\kappa N}{2}\right)$; the noise slows the wave speed. If the stochastic integral is taken instead in the sense of Stratonovich, then the asymptotic wave speed is the classical McKean wave speed and does not depend on $\epsilon$. 

In the second situation (noise with integrable spatial covariance, stochastic integral taken in the sense of Itô), a travelling front can be defined for all $\epsilon > 0$. Its average asymptotic speed does not depend on $\epsilon$ and is the classical wave speed of the unperturbed KPP equation. 

 \end{abstract}
 

\section{Introduction} This article considers the Kolmogorov Petrovskii Piscunov equation with stochastic perturbation:

\begin{equation}\label{eqKPP}
 \left\{\begin{array}{lll} \partial_t u &=&
\left(\frac{\kappa}{2}u_{xx} + u(1 - u)\right) dt +   \epsilon u 
\partial_t \zeta 
\\ u_0 (x) &=& {\bf 1}_{(-\infty, -\frac{1}{N}\log 2)}(x) + \frac{1}{2}e^{-Nx}
{\bf 1}_{[-\frac{1}{N}\log 2, +\infty)} (x) \end{array}\right.
\end{equation}

\noindent where $\zeta$ is a Gaussian random field, adapted to a filtered probability space $(\Omega, {\cal G}_t, {\cal G}, \mathbb{Q})$, satisfying: \[ \mathbb{E}_{\mathbb{Q}} \left [ \zeta \right ] =  0, \qquad \mathbb{E}_{\mathbb{Q}} \left [ \zeta(s,x)\zeta(t,y) \right ] = (s \wedge t) \Gamma (x-y)\]  where $\Gamma  \in C^\infty (\mathbb{R}^d)$. Two situations are considered; firstly, $\Gamma \equiv 1$ and secondly $\int_{\mathbb{R}^d} |\Gamma(z)|dz < +\infty$.  

\paragraph{Background} When $\epsilon = 0$, the equation under discussion is the travelling wave of the Kolmogorov
Petrovskii Piskunov equation

\[ u_t = \frac{\kappa}{2} u_{xx} + u(1-u)\]

\noindent introduced in 1937 by both Fisher~\cite{Fish} and Kolmogorov, Petrovskii and Piscunov~\cite{KPP}. It has received
substantial attention since then; a probabilistic interpretation of the solution in terms of branching Brownian motion was given by McKean~\cite{McK}(1975) and~\cite{McK2}(1976). Refinements and more precise description of the asymptotics were given by Bramson~\cite{Br78}(1978) and developed further by Bramson~\cite{Br}(1983).   For the KPP equation (when $\epsilon = 0$ in Equation~\eqref{eqKPP}), there is a one parameter family $F_\gamma : \gamma \geq \sqrt{2\kappa}$ of  travelling front solutions $F_\gamma (x - \gamma t)$ where $F(x) \simeq e^{-\nu x}$ for large $x$ and $\gamma = \frac{\kappa \nu}{2} + \frac{1}{\nu}$, $\nu \leq \sqrt{\frac{2}{\kappa}}$. For initial condition $u_0$ such that $u_0(x) \leq 1$ for $x \leq 0$ and $u_0(x) \leq \exp\left\{-\sqrt{\frac{2}{\kappa}}x \right\}$ there is convergence to the travelling front with minimal speed $\gamma_0 = \sqrt{2\kappa}$, in the sense that 

\[\sup_{x \in \mathbb{R}}\left(\lim_{t \rightarrow +\infty} \left | u(t,x) - F_{\gamma_0}(x - g(t)) \right |\right) = 0\]

\noindent where $g(t)$ is the {\em marker}, defined by: $u(t,g(t)) = \frac{1}{2}$, and $F_{\gamma_0}(0) = \frac{1}{2}$. The marker satisfies: 

\[ \lim_{t \rightarrow +\infty} \frac{g(t)}{t} = \gamma_0.\] 

\paragraph{Results} The results of this article are stated in the abstract; with  $\Gamma(x) \equiv 1$ (so that $\zeta$ is simply a standard Wiener process, the same Wiener process for each $x \in \mathbb{R}$), there is slow-down of the travelling front; the asymptotic speed is $\max\left(\sqrt{2\kappa \left (1 - \frac{\epsilon^2}{2} \right )}, \frac{1}{N}(1 - \frac{\epsilon^2}{2})+ \frac{\kappa N}{2}\right)$ and the front breaks down (because the solution tends to zero) if $\epsilon > \sqrt{2}$. If the stochastic term is taken in the sense of Stratonovich, however, there is no slow-down and the average wave speed is the same for all $\epsilon > 0$. This could naïvely suggest that the slow-down is simply a consequence of taking a wrong approach to the way that the noice is introduced and that a Stratonovich integral is correct. 

This naïve suggestion is, however, flatly contradicted by the results for   $\Gamma$ satisfying $\int |\Gamma(x)|dx < +\infty$, with the stochastic integral taken in the sense of Itô. In this setting, the average speed of the travelling front does not change with $\epsilon$.  

The results for the speed of the travelling front, when $\Gamma \equiv 1$, are already known from the work of Øksendal, Våge and Zhao (2001)~\cite{OkVaZh}, although the treatment presented here is different and uses a different definition of the wave marker, which leads to slightly sharper results. The outstanding question of interest dealt with in this article, therefore, is how a decay  in the spatial covariance function for the noise affects the results. The answer is that the spatial variation in the noise introduces a drift which precisely compensates the Øksendal-Våge-Zhao decay (although I don't have any intuitive explanation for why they should exactly cancel), so that the asymptotic speed of the travelling front is the same as that for the unperturbed equation, when the stochastic integral is taken in the sense of Itô.

 \paragraph{Motivation} My personal motivation for studying this came from the results of Mueller, Mytnik, Quastel~\cite{MuMyQu} where the equation

\begin{equation} \label{eqmumyqu} \left\{\begin{array}{l} u_t = \frac{1}{2}u_{xx} + u(1-u) + \epsilon u^{1/2} \xi \\ u(0,x) = {\bf 1}_{(-\infty, 0]} \end{array} \right. \end{equation}

\noindent was considered, where $\xi$ is space time white noise and the stochastic integral is taken in the sense of Walsh~\cite{Wa}. The `noise' term in Mueller et. al. is entirely different from that under consideration here. The techniques of proof are entirely different, since with space/time white noise, spatial derivatives do not exist (I make strong use of them in this article) and the Wiener sheet does not have sufficient regularity to allow for an Itô formula (I make heavy use of Itô's formula here). The work of Mueller, Mytnik and Quastel requires entirely different and much more sophisticated mathematical techniques for proving the results and the results are entirely different.  

For the noise in Mueller et. al., the slow-down is {\em greater} than that of Øksendal-Våge-Zhao; the slow-down follows the Derrida-Brunet conjecture. That is, if $\gamma_0$ denotes the speed of the unperturbed wave front and $\gamma_\epsilon$ the speed of the wave front with perturbation parameter $\epsilon$, then

\[ \gamma_0 - \gamma_\epsilon \sim \frac{\pi^2}{4 \left(\log\frac{1}{\epsilon}\right)^2}\]

\noindent as $\epsilon \rightarrow 0$. 

The `noise' term differs in two respects: firstly, it  has $ u^{1/2}$ instead of $u$. This is already a major difference; parabolic SPDEs with $u^{1/2}\xi$ noise are known to have markedly different properties than SPDEs with $u \xi$ noise.   Also, the noise is a space/time white noise instead of having smooth space correlation. The behaviour of SPDEs of this type, whose initial conditions have compact support is markedly different; it is well known that the solutions themselves have compact support for all time. Indeed, the same carries over when there is a right limit for the support of a bounded initial condition and the definition of the marker in~\cite{MuMyQu} is the right end point of the support. 

This article shows that the slow-down is heavily dependent on the structure of the noise; this article, taken together with~\cite{MuMyQu}, illustrates that the features of the noise responsible for the slow-down are not well understood at all; this is an important outstanding problem. This article makes no attempt to analyse why the noise in these three situations (the two presented here and the situation in~\cite{MuMyQu}) give such markedly different results. This question is left for further study.

\section{Standard Wiener Process Noise}\label{secseason}

In this section, Equation~\eqref{eqKPP} is considered where the noise  is simply a standard Wiener process; $\Gamma \equiv 1$. This equation has been considered by Øksendal et. al.~\cite{OkVaZh}. In that paper, a marker of the form $g(t) = \gamma t$ is considered. The correct value of $\gamma$ is obtained, although the description of the wave front is not sharp. The approach here defines the marker differently and presents an approach which gives a sharper description of the travelling wave and is hopefully more transparent.  Let $W$ denote a standard Wiener process. Let $u$ and $v$ solve:

\begin{equation}\label{eqseason} \left\{ \begin{array}{l} \partial_t u = \left(\frac{\kappa}{2} u_{xx} + u(1-u)\right)dt + \epsilon u dW \\ u(0,x) = {\bf 1}_{(-\infty, -\frac{1}{N}\log 2)(x)} + e^{-Nx}{\bf 1}_{[-\frac{1}{N}\log 2, +\infty)}(x) =:u_N \\ d v = v(1-v)dt + \epsilon v d W \\ v(0) \equiv 1 \end{array}\right.
\end{equation}

\noindent Throughout, the notation $\partial_t$ will denote a differential with respect to the variable denoted $t$ when there are several variables; $d$ denotes the differential when there is no spatial dependence. The stochastic integral is taken in the sense of Itô.  

\begin{Lmm}\label{lmmvinv}  Consider $v$ from~\eqref{eqseason}.  Provided $0 < \epsilon < \sqrt{2}$, $V:= \lim_{t \rightarrow +\infty} \mathbb{E}[v(t)] = 1 - \frac{\epsilon^2}{2}$  and $v$  has a non-trivial stationary distribution where the Laplace functional of the limiting distribution is given by:

 \[ \mathbb{E} \left [e^{\lambda v}\right] = \left(1 - \frac{\epsilon^2}{2}\lambda \right)^{1 - (2/\epsilon^2)}.\]
 
\noindent If $\epsilon \geq \sqrt{2}$, then $\lim_{t \rightarrow +\infty}  v(t)  = 0$ almost surely.  
\end{Lmm}

\paragraph{Sketch of Proof} Only a sketch is given, since it is routine and the result is well known. For the Itô formulation, standard Itô calculus gives:

\[ \frac{\partial}{\partial t} \mathbb{E} \left [ e^{\lambda v(t)}\right ] = \lambda   \mathbb{E} \left [ e^{\lambda v(t)}v(t)\right ] - \lambda \mathbb{E} \left [ e^{\lambda v} v^2 \right ] + \frac{\epsilon^2 \lambda^2}{2} \mathbb{E} \left [ e^{\lambda v} v^2 \right].\]

\noindent Let $M(t,\lambda) = \mathbb{E} \left [ e^{\lambda v(t)}\right ]$, then $M$ satisfies:

\[ \frac{\partial}{\partial t} M(t,\lambda) =  \lambda \frac{\partial}{\partial \lambda} M(t,\lambda) - \lambda \frac{\partial^2}{\partial \lambda^2} M(t,\lambda) + \frac{\epsilon^2 \lambda^2}{2} \frac{\partial^2}{\partial \lambda^2} M(t,\lambda).\]

\noindent Let $M(\lambda) = \lim_{t \rightarrow +\infty} M(t,\lambda)$, then $M(\lambda)$ satisfies:

\[ 0 =  \frac{\partial}{\partial \lambda} M(\lambda) - \frac{\partial^2}{\partial \lambda^2} M(\lambda) + \frac{\epsilon^2 \lambda}{2} \frac{\partial^2}{\partial \lambda^2} M(\lambda).\] 

 \noindent Let $f(\lambda) = M^\prime(\lambda)$, then $f$ satisfies:
 
 \[ \frac{\partial}{\partial \lambda} \log f(\lambda) = \frac{2}{2 - \epsilon^2 \lambda} \]
 
 \noindent giving 
 
 \[ f(\lambda) = f(0) \left(\frac{2}{2 - \epsilon^2 \lambda}\right)^{2/\epsilon^2} \qquad \lambda < \frac{2}{\epsilon^2}.\]
 
 \noindent Using $M(0) = 1$ and $V := \lim_{t \rightarrow +\infty}\mathbb{E}[v(t)] = f(0)$, it follows that:

\[ M(\lambda) = 1 + \frac{V}{1 - \frac{\epsilon^2}{2}} \left(   \frac{1}{(1 - \frac{\epsilon^2 \lambda}{2})^{(2 /\epsilon^2) - 1}} - 1\right). \]
 
 \noindent Either $v(t) \rightarrow 0$ in law (in which case $\lim_{\lambda \rightarrow -\infty} \mathbb{E} \left [ e^{\lambda v} \right ] = 1$), or else $\lim_{\lambda \rightarrow -\infty} \mathbb{E} \left [e^{\lambda v} \right ] = 0$, from which it follows that for $\epsilon < \sqrt{2}$, $V = 1 - \frac{\epsilon^2}{2}$ and  the result follows. \qed \vspace{5mm}

\noindent In this section, with the noise taken as a standard Wiener process, heavy use will be made of the {\em normalised} solution, defined as $\widetilde{u} := \frac{u}{v}$, 

\begin{Lmm} Let $(u,v)$ satisfy Equation~\eqref{eqseason}. Let $\widetilde{u}:= \frac{u}{v}$, then $\widetilde{u}$ satisfies:

\begin{equation}\label{eqtildu} \left\{ \begin{array}{l} \widetilde{u}_t = \frac{\kappa}{2} \widetilde{u}_{xx} + v(t) \widetilde{u} (1 - \widetilde{u}) \\ \widetilde{u}(0,x) = {\bf 1}_{(-\infty, -\frac{1}{N}\log 2)}(x) + \frac{1}{2}e^{-Nx}
{\bf 1}_{[-\frac{1}{N}\log 2, +\infty)} (x) \end{array}\right.
\end{equation}
\end{Lmm}

\paragraph{Proof} Straightforward application of Itô's formula. \qed \vspace{5mm}

\begin{Lmm}\label{lmmnonincrease} The function $\widetilde{u}(t,.)$, defined by Equation~\eqref{eqtildu} is a strictly decreasing  function for each $t \geq 0$.
\end{Lmm}

\paragraph{Proof} Set $\widetilde{\theta} = -(\log \widetilde{u})_x$. Then $\widetilde{\theta}$ satisfies:

\begin{equation}\label{eqtilth} \left\{ \begin{array}{l} \widetilde{\theta}_t = \frac{\kappa}{2} \widetilde{\theta}_{xx} - \kappa \widetilde{\theta} \widetilde{\theta}_x - v\widetilde{u} \widetilde{\theta} \\ \widetilde{\theta}(0,x) = N {\bf 1}_{[-\frac{1}{N} \log 2, +\infty)}(x). \end{array}\right.
\end{equation}

\noindent Since the initial condition is non-negative and equal to $N$ for $x < -\frac{1}{N} \log 2$, it is straightforward that solutions are strictly positive for all $t > 0$. \qed\vspace{5mm}

\noindent For the problem~\eqref{eqseason}, the marker is defined as follows:

\begin{Defn}[$a$-Marker]\label{defsmark} Let $\widetilde{u}$ satisfy Equation~\eqref{eqtildu}. For $a \in (0,1)$, the {\em  $a$-marker} is defined as the random, adapted function $g^{(a)}$ that satisfies: 

\[ \widetilde{u}(t, g^{(a)}(t)) \equiv a \qquad \forall t > 0. \]
\end{Defn}

\begin{Lmm}\label{Lmmamark} For each $a \in (0,1)$ the $a$-marker is well defined.
\end{Lmm}

\paragraph{Proof of Lemma~\ref{Lmmamark}} This follows directly from Lemma~\ref{lmmnonincrease}. \qed\vspace{5mm} 

\noindent {\bf Note} From Equation~\eqref{eqtildu}, it follows that

\begin{equation}\label{eqttil2}
0 = \frac{\kappa}{2} \widetilde{u}_{xx}(t,g^{(a)}(t)) + \dot{g}^{(a)}\widetilde{u}_x(t,g^{(a)}(t) + v(t)a(1-a) 
\end{equation}

\noindent so that $\dot{g}^{(a)}$ satisfies:

\[ \dot{g}^{(a)}(t) = \frac{\frac{\kappa}{2}\widetilde{u}_{xx}(t,g^{(a)}(t)) + v(t)a(1-a)}{-\widetilde{u}_x(t,g^{(a)}(t))}. \]

\begin{Lmm} Let $\widetilde{u}$ satisfy Equation~\eqref{eqtildu} and $g^{(a)}$ follow  Definition~\ref{defsmark}. Let $\widetilde{\theta} = -(\log \widetilde{u})_x$ so that $\widetilde{\theta}$ satisfies Equation~\eqref{eqtilth}. Then:
\begin{enumerate} \item  $\lim_{t \rightarrow +\infty} \frac{g^{(a)}(t)}{t} =: \overline{g}^{(a)}$ exists, in the sense that there is a number $\overline{g}^{(a)} \in \mathbb{R}_+ \cup \{+\infty\}$ such that $\lim_{t \rightarrow +\infty} \frac{g^{(a)}(t)}{t} = \overline{g}^{(a)}$ almost surely. 
\item There is a value $\widetilde{\theta}_N^{(a)} \in \mathbb{R}$ such that 

\begin{equation}\label{eqthlim} \lim_{x \rightarrow +\infty} \left(\liminf_{t \rightarrow +\infty}\widetilde{\theta}(t, g^{(a)}(t) + x)\right) = \lim_{x \rightarrow +\infty}\left(\limsup_{t \rightarrow +\infty} \widetilde{\theta}(t,g^{(a)}(t) + x)\right) =  \widetilde{\theta}_N^{(a)}
\end{equation}
almost surely. 
\end{enumerate}
\end{Lmm}

\paragraph{Proof} The proof of these two points requires ergodic theorems that are found in Chapter 20 of Kallenberg~\cite{Kal}. Let   $(v^{(T)} ,\widehat{u}^{(T,a)})$ be the solution to:

\begin{equation}\label{eqbackeq}
\left\{ \begin{array}{l} \widehat{u}_t^{(T,a)} = \frac{\kappa}{2}\widehat{u}_{xx}^{(T,a)} +  \widehat{u}^{(T,a)}_x\dot{g}^{(T,a)} + v^{(T)}(t)\widehat{u}^{(T,a)}(1 - \widehat{u}^{(T,a)}) \\ 
\dot{g}^{(T,a)}(t) = - \frac{a(1-a)v(t)}{\widehat{u}_x^{(T,a)}(t,0)} - \frac{\frac{\kappa}{2} \widehat{u}_{xx}^{(T,a)}(t,0)}{\widehat{u}_x^{(T)}(t,0)}\\
dv^{(T)} = v^{(T)}(1-v^{(T)})dt + \epsilon v^{(T)} dW\\
\widehat{u}^{(T,a)}(-T,x) = {\bf 1}_{(-\infty, -\frac{1}{N}\log \frac{1}{a})}(x) + a e^{-Nx} {\bf 1}_{[-\frac{1}{N}\log \frac{1}{a}, +\infty)}(x)\\ 
v^{(T)} (-T) = 1
\end{array} \right. 
\end{equation}

\noindent where $W$ is a Wiener process started at time $-T$.

\paragraph{Note} From the definition, it follows that 

\[\widehat{u}^{(T,a)}(t,x) = \widehat{u}^{(T,1/2)}(t,x + u^{(T,1/2)-1}(t,a)) \qquad \forall a \in (0,1) \]

\noindent where 

\[ \widehat{u}^{(T,a)-1}(t,x) := \{y : \widehat{u}^{(T,a)}(t,y) = x\}.\] 

\noindent Let $\nu^{(T)}$ denote the probability distribution of $(v^{(T)}(0), u^{(T)}(0,.))$ and let $\nu_N$ denote any probability distribution which is a limit point of $(\nu^{(T)})_{T > 0}$. Existence of such a limit may be established using Prohorov's theorem and the Prohorov metric. Let 

\[{\cal L} = \left\{ u : u \quad \mbox{non-increasing}, \quad u(x) \stackrel{x \rightarrow +\infty}{\longrightarrow} 0, \quad u(x) \stackrel{x \rightarrow -\infty}{\longrightarrow} 1\right \}.\] Now let $S = \mathbb{R}_+ \cup \{+\infty\} \times {\cal L}$ and let $d$ denote the metric;

\[ d((v_1,u_1),(v_2,u_2)) = \sqrt{ \left (\frac{1}{1+v_1} - \frac{1}{1 + v_2} \right )^2 + \frac{1}{2}\int_{-\infty}^\infty e^{-|x|} (u_1(x) - u_2(x))^2 dx} \]

\noindent It is straightforward that $(S,d)$ is compact; any Cauchy sequence of functions $u \in {\cal L}$ clearly has a limit and, for the first argument, $d$ was chosen to make $\mathbb{R}_+ \cup \{+\infty\}$ compact. Now let $({\cal P}(S,d),\rho)$ denote the space of probability measures over $(S,d)$ endowed with Prohorov metric $\rho$, which is defined as follows. Let ${\cal S}$ denote the Borel sets of $S$ and let 

\[ \rho(\lambda_1, \lambda_2) = \inf\{\alpha > 0 : \lambda_1 (A) \leq \lambda_2(A_\alpha) + \alpha, \quad \lambda_2(A) < \lambda_1(A_\alpha) + \alpha \qquad \forall A \in {\cal S}\} \]

\noindent where 

\[A_\alpha = \{(v,u(.))\in \mathbb{R}_+ \times {\cal L} : \inf_{(a,b(.)) \in A} d((v,u(.)),(a,b(.))) \leq \alpha\}.\]

\noindent  It is a standard result that if a space $(S,d)$ is compact, then the space of probability measures over $(S,d)$ endowed with Prohorov metric $({\cal P}(S,d), \rho)$ is compact. It follows that the sequence $(\nu^{(T)})_{T > 0}$ has a converging subsequence with limit $\nu_N$, which is a probability measure over $(S,d)$.  

Having established existence of a limit $\nu_N$, it follows directly from the fact that $(v,\widetilde{u}(.))$ is a Feller process that this limit is unique and that $\rho(P_t((1,u_N(.)),.),\nu_N(.))\stackrel{t \rightarrow +\infty}{\longrightarrow} 0$ where $P_t(x,.)$ denotes the transition semigroup of the Feller process $(v(t), \widehat{u}(t,.))_{t \geq 0}$. Let 

\[ S_N = \overline{\cup_{t \geq 0} \mbox{suppt}(P_t((1,u_N(.)),.)}\]

\noindent namely, the closure (under metric $d$) of the union of the supports of the probability measures \\ $P_t((1,u_N(.)),.)$. Let ${\cal S}_N$ denote the Borel $\sigma$-algebra of $S_N$. \vspace{5mm}

\noindent Following Kallenberg~\cite{Kal}, page 405 Theorem 20.17, a regular Feller process is {\em either} Harris recurrent, {\em or} uniformly transient. The definition of uniformly transient is given just above the statement of Theorem 20.17 in Kallenberg~\cite{Kal}; the processs $(v(t),\widehat{u}(t,.))$ is uniformly transient if 

\[ \sup_{x \in S_N} \mathbb{E}_x \left [ \int_0^\infty {\bf 1}_K ((v(t),\widehat{u}(t,.))) dt \right ] < +\infty \qquad \forall K \subset {\cal S}_N \quad \mbox{compact}.\]

\noindent This clearly does not hold since for any compact $K$ with $\nu_N(K) > 0$, 
\[ \lim_{t \rightarrow +\infty} \frac{1}{t} \sup_{x \in S_N} \mathbb{E}_x \left [ \int_0^t {\bf 1}_K((v(s),\widehat{u}(s,.))ds \right ] = \nu_N(K) > 0.\]

\noindent It follows that $(v(t), \widehat{u}(t,.))$ is Harris recurrent with supporting measure $\nu_N$ (the definition is given by Equation (6) on page 400 of~\cite{Kal}). Hence, by Theorem 20.12 page 400 of Kallenberg~\cite{Kal}, it is strongly ergodic. 
\vspace{5mm}

\noindent It follows that $(v(t), \widehat{u}(t,.))_{t \geq 0}$ is a Feller process on $(S,d)$ which, by Theorem 20.21 page 409 of Kallenberg~\cite{Kal} is positive recurrent with invariant distribution $\nu_N$ and, for any functional $F$,

\begin{equation} \label{eqerglim}\lim_{t \rightarrow +\infty} \frac{1}{t} \int_0^t F((v(s),\widehat{u}(s,.))ds = \mathbb{E}_{\nu_N} \left [ F((v,\widehat{u}(.))\right ].
\end{equation}

\paragraph{Showing that $\overline{g}^{(a)}$ is well defined in $\mathbb{R}_+ \cup \{+\infty\}$} To show that $\overline{g}^{(a)}$ is well defined, consider the functional:

\[ F(v,\widehat{u}) = - \frac{a(1-a)v}{\widehat{u}_x(0)} - \frac{\frac{\kappa}{2} \widehat{u}_{xx}(0)}{\widehat{u}_x(0)}.\]  

\noindent The result now follows directly from Equation~\eqref{eqerglim}, because from Equation~\eqref{eqbackeq}, it follows that 

\begin{equation}\label{eqgdot} \dot{g}^{(a)}(t) = - \frac{a(1-a)v(t)}{\widehat{u}_x(t,0)} - \frac{\frac{\kappa}{2} \widehat{u}_{xx}(t,0)}{\widehat{u}_x(t,0)}.
\end{equation}

\paragraph{Showing that $\lim_{x \rightarrow +\infty} \left(\lim_{t \rightarrow +\infty} (-\log \widehat{u}(t,x))_x\right)$ is well defined}

Let $(v, \widehat{u})$ denote the solution to Equation~\eqref{eqbackeq} with $T = 0$. Let $\widetilde{\theta}(t,x) = -(\log \widehat{u}(t,x))_x$. It follows from the ergodic theorem that

\[ \lim_{t \rightarrow +\infty}\frac{1}{t} \int_0^t \widetilde{\theta}(s,x)ds = \mathbb{E}_{\nu_N}[\widetilde{\theta}(x)], \]

\noindent where $\widetilde{\theta}(x)$ (without a $t$ dependence) denotes: $\widetilde{\theta}(x) = -(\log \widehat{u})_x$, where $(v,\widehat{u}) \sim \nu_N$.\vspace{5mm}

\noindent Let $\widetilde{\theta}_N = \lim_{x \rightarrow +\infty} \mathbb{E}_{\nu_N} [\widetilde{\theta}(x)]$. Now consider $\theta$ satisfying

\[ \left\{ \begin{array}{l} dv_t = v(1-v) dt + \epsilon v dW \\ \theta_t = \frac{\kappa}{2} \theta_{xx} - \kappa \theta \theta_x + \dot{g} \theta_x - v\widehat{u} \theta \\
            \theta_0 = -(\log \widehat{u}_0)_x, \quad (v(0),\widehat{u}_0) \sim \nu_N
           \end{array} \right. \]

\noindent Let $B$ denote a Wiener process, independent of $W$, with diffusion coefficient $\kappa$; let $({\cal X}, {\cal F}, ({\cal F}_t)_{t \geq 0}, \mathbb{P})$ denote the probability space for $B$ and let $\mathbb{E}_{\mathbb{P}}$ denote expectataion with respect to $\mathbb{P}$. Then $\theta$ has representation:

\begin{equation}\label{eqthetarep} \theta(t,x) = \mathbb{E}_{\mathbb{P}}\left [ \theta(0,X_{0,t}(x))\exp\left\{ -\int_0^t v(s)\widehat{u}(s,X_{s,t}(x))ds\right\} \right ]
\end{equation}

\noindent where:

\[ X_{s,t}(x) = x + (B(t) - B(s)) + (g(t)-g(s)) - \kappa \int_s^t \theta(r,X_{r,t}(x))dr.\]

\noindent Let $\theta^* = \lim_{x \rightarrow +\infty} \theta(0,x)$. It follows directly from the formula~\eqref{eqthetarep} that $\lim_{x \rightarrow +\infty} \sup_t \theta(t,x) \leq \theta^*$. 

Now suppose that the initial condition $(v_0, \widehat{u})  \in \mbox{suppt}(\nu_N)$ satisfies: $\lim_{x \rightarrow +\infty} -(\log \widehat{u})_x(x) = \theta^* < \widetilde{\theta}_N$ where the inequality is strict. It follows that 

\[ \lim_{x \rightarrow +\infty} \left(\lim_{t \rightarrow +\infty}\frac{1}{t}\int_0^t \theta(s,x)ds\right) \leq \theta^* < \widetilde{\theta}_N\]

\noindent which contradicts the ergodic theorem. It follows that, for $(v,\widehat{u}) \in \mbox{suppt}(\nu_N)$ and $\theta = -(\log \widehat{u})_x$, $\lim_{x \rightarrow +\infty} \mathbb{E}_{\nu_N} [ \theta(x) ] = \widetilde{\theta}_N$ while almost surely, $\lim_{x \rightarrow +\infty} \theta(x) \geq \widetilde{\theta}_N$. By Fatou's lemma, it now follows that 

\[ \widetilde{\theta}_N = \lim_{x \rightarrow +\infty} \mathbb{E}_{\nu_N}[\theta(x)] \geq \mathbb{E} [\lim_{x \rightarrow +\infty} \theta(x)] \geq \widetilde{\theta}_N \]

\noindent and hence for $(v, \widehat{u}) \in \mbox{suppt}(\nu_N)$, $\lim_{x \rightarrow +\infty} -(\log \widehat{u})_x(x) = \widetilde{\theta}_N$, $\nu_N$ almost surely and the result follows.\qed \vspace{5mm}

\noindent Let $\widetilde{u}$ denote the solution to Equation~\eqref{eqtildu} and let $\widetilde{\theta} = -(\log \widetilde{u})_x$. Then $\widetilde{u}_{xx} = \widetilde{u}(\widetilde{\theta}_x - \widetilde{\theta}^2)$ and $\widetilde{u}_x = - \widetilde{u} \widetilde{\theta}$.  Using this, it follows directly from~\eqref{eqttil2} that $\dot{g}^{(a)}$ satisfies:  

\begin{equation}\label{eqdotg} \dot{g}^{(a)}(t) = -\frac{\kappa}{2} (\log \widetilde{\theta})_x(t,g^{(a)}(t)) + \frac{\kappa}{2}\widetilde{\theta} (t,g^{(a)}(t)) + \frac{v(t)(1-a)}{\widetilde{\theta} (t,g^{(a)}(t))}.
\end{equation}

\begin{Lmm} \label{lmmmark} The wave marker satisfies:
\begin{equation}
\overline{g} \leq \left\{ \begin{array}{ll} \sqrt{2\kappa \overline{v}} & N > \sqrt{\frac{2\overline{v}}{\kappa}} \\ \frac{\kappa N}{2} + \frac{\overline{v}}{N} & N \leq \sqrt{\frac{2\overline{v}}{\kappa}}\end{array}\right. 
\end{equation}
\end{Lmm}

\paragraph{Proof} Let $\widetilde{u}(t,x)$ satisfy Equation~\eqref{eqtildu} and let   $g(t):= g^{(1/2)}(t)$ (the marker with $a = \frac{1}{2}$). Then $U(t,x):= \widetilde{u}(t,g(t) + x)$ has a Kacs-Feynman representation:  

\[ U(t,x) = \mathbb{E}_{\mathbb{P}} \left [\widetilde{u}(0,x + B(t) + g(t)) e^{\int_0^t v(s)( 1 - \widehat{u}(s, x+g(t) + B(t) - B(s)))ds} \right ] \]

\noindent where $\mathbb{P}$ is a probability measure under which $B$ is a Brownian motion with diffusion coefficient $\kappa$.  An upper bound may be obtained quite easily using the fact that $B(t) \sim N(0, \kappa t)$;

\begin{eqnarray*}\lefteqn{U(t,x) \leq e^{\int_0^t v(s)ds}}\\&& \times \left(\int_{-\infty}^{-\frac{1}{N}\log 2 -x - g(t)} \frac{1}{\sqrt{2\pi \kappa t}}e^{-y^2/2\kappa t} dy + \int_{-\frac{1}{N}\log 2 - x - g(t)}^\infty \frac{1}{\sqrt{2\pi \kappa t}}e^{-N(x + g(t) + y) - \frac{y^2}{2\kappa t}} dy \right).
\end{eqnarray*}

\noindent Let $\Phi(z) = \mathbb{P}(Z \leq z)$ where $Z \sim N(0,1)$. Then:

\begin{eqnarray}\label{equtildub} \lefteqn{U(t,x) \leq e^{\int_0^t v(s)ds}}\\&& \times \nonumber \left( \Phi\left(-\frac{g(t)}{\sqrt{\kappa t}} - \frac{x + \frac{1}{N} \log 2}{\sqrt{\kappa t}}\right) + e^{-Nx + \frac{\kappa t N^2}{2} - Ng(t)} \Phi\left (\frac{g(t)}{\sqrt{\kappa t}} + \sqrt{\kappa t} N + \frac{\frac{1}{N} \log 2 + x}{\sqrt{\kappa t}}   \right ) \right ) \end{eqnarray}

\noindent By construction, $U(t,0) \equiv \frac{1}{2}$ and $U(t,x) > \frac{1}{2}$ for $x < 0$. These provide lower bounds for the right hand side of Equation~\eqref{equtildub}. 

\begin{itemize}
\item The first term tends to zero if   $\overline{g} > \sqrt{2\kappa \overline{v}}$.
\item The second term tends to zero if $\sqrt{2\kappa \overline{v}} < \overline{g} < \kappa N$ and $\left\{\kappa N < \overline{g} \right \} \cap \left\{ \overline{g} > \frac{\kappa N}{2} + \frac{\overline{v}}{N}\right\}$.  
\end{itemize}

\noindent It follows that, for the condition $U(t,0) \equiv \frac{1}{2}$ and $U(t,x) > \frac{1}{2}$ for all $x < 0$ to hold, it is necessary that

\[ \overline{g} \leq  \left\{ \begin{array}{ll} \sqrt{2\kappa \overline{v}} & N \geq \sqrt{\frac{2\overline{v}}{\kappa}} \\ \frac{\kappa N}{2} + \frac{\overline{v}}{N} & N \leq \sqrt{\frac{2\overline{v}}{\kappa}} \end{array}\right. \]

\noindent and Lemma~\ref{lmmmark} follows. \qed

\begin{Lmm} 

\begin{enumerate} 
\item Provided $1 > \frac{\epsilon^2}{2}$ in~\eqref{eqseason}, it follows that 
 For all $a,b$ such that $0 < a < b < 1$, 
 \begin{equation}\label{eqgagb}\lim_{N \rightarrow +\infty}  \left(\lim_{t \rightarrow +\infty}  \mathbb{P}\left(\sup_{0 \leq s \leq t} \left | g^{(a)}(s)  -  g^{(b)}(s) \right | > N \right) \right ) = 0.\end{equation}
 
 \noindent From this, it follows that $\overline{g}^{(a)} = \overline{g}$ (the same value) for each $a \in (0,1)$ and $\widetilde{\theta}_N^{(a)} = \widetilde{\theta}_N$ (the same value) for each $a \in (0,1)$.

\item In Equation~\eqref{eqdotg}, $\lim_{a \downarrow 0} \left(\lim_{t \rightarrow +\infty} (\log \widetilde{\theta})_x(t,g^{(a)}(t))\right) = 0$ almost surely. 
\item  $\overline{g}$  satisfies:

\begin{equation}\label{eqovg} \overline{g} = \frac{\kappa}{2}\widetilde{\theta}_N + \frac{\overline{v}}{\widetilde{\theta}_N}
\end{equation}

\noindent where $\overline{v} = \lim_{t \rightarrow +\infty} \frac{1}{t}\int_0^t v(s) ds = \lim_{t \rightarrow +\infty} \mathbb{E}[v(t)] = 1 - \frac{\epsilon^2}{2}$.
\end{enumerate}
\end{Lmm}

\paragraph{Proof} 
\begin{enumerate}
 \item   Suppose there exists an  $a$ and a $b$ such that $a < b$ (so that $g^{(a)} \geq g^{(b)}$) and a sequence of times $t_n \rightarrow +\infty$ such that $g^{(a)}(t_n) - g^{(b)}(t_n) \stackrel{n \rightarrow +\infty}{\longrightarrow} +\infty$. Now recall $\widetilde{u}^{(a)}(t,.):= \widetilde{u}(t,g^{(a)}(t) + .)$. Then $\limsup_{x \rightarrow -\infty}\left(\limsup_{n \rightarrow +\infty} \widetilde{u}^{(a)}(t_n,x)\right) \leq b < 1$, while $\lim_{n \rightarrow +\infty}\widetilde{u}^{(a)}(t_n,0) = a$ and, since $\widetilde{u}^{(a)}(t,.)$ is decreasing in $x$ for each $t > 0$, it follows that 
 
 \[ 0 < a \leq \liminf_{x \rightarrow -\infty}\left(\liminf_{n \rightarrow +\infty} \widetilde{u}^{(a)}(t_n,x)\right) \leq b.\]

 \noindent Let ${\cal L}^{(a,T)}$ of the pair $((\widetilde{u}^{(a,T)}(t,.),v(t)): t \geq 0)$ where $(\widetilde{u}^{(a,T)},v)$ satisfy:
 
 \[ \left\{ \begin{array}{l} \widetilde{u}_t^{(a,T)} = \frac{\kappa}{2}\widetilde{u}_{xx}^{(a,T)} + \dot{g}^{(a)}\widetilde{u}_x^{(a,T)} + v \widetilde{u}^{(a,T)}(1 - \widetilde{u}^{(a,T)}) \\ \widetilde{u}^{(a,T)}(-T,.) = {\bf 1}_{(-\infty, -\frac{1}{N}\log\frac{1}{a}]} + e^{-Nx}{\bf 1}_{(-\frac{1}{N}\log\frac{1}{a},+\infty)} \\ \partial_t v = v(1-v)dt + \epsilon v dW_t \\ v(-T) = 1 \end{array}\right.  \]

 \noindent (note that we are considering $(\widetilde{u}^{(a,T)}(t,.),v(t)):  t \geq 0$; the solutions from time $t = 0$ onwards). Let ${\cal L}^{(a)}$ denote the limit of ${\cal L}^{(a,T)}$. Now let $A_\epsilon = \{ (\widetilde{u}^{(a)},v) : \lim_{t \rightarrow +\infty} \inf_{s \geq t}\lim_{x \rightarrow -\infty} \widetilde{u}^{(a)}(s,x) \leq 1 - \epsilon\}$. Let $\phi(s) = \lim_{x \rightarrow -\infty}\widetilde{u}^{(a)}(s,x)$ and note that $\phi(s) \geq a$. Furthermore, $\phi(s)$ solves $\dot{\phi}(s) = v(s)\phi(s)(1-\phi(s))$ and hence (using $\phi(s) \geq a > 0$) it follows that either $\phi(s) \uparrow 1$ as $s \rightarrow +\infty$ or $v(s) \stackrel{s \rightarrow +\infty}{\longrightarrow} 0$.  Since this event has probability $0$, it follows that ${\cal L}^{(a)}(A_\epsilon) = 0$ for all $\epsilon > 0$ and hence~\eqref{eqgagb} follows.   \vspace{5mm}

\noindent From this, it is immediate that  $\overline{g}^{(a)}$ takes the same value for each $a \in (0,1)$. 

It now follows directly from the fact that 

\[ \widetilde{\theta}_N^{(a)} = \lim_{x \rightarrow +\infty} \left(\liminf_{t \rightarrow +\infty} \widetilde{\theta}(t, g^{(a)}(t) + x)\right) = \lim_{x \rightarrow +\infty} \left(\limsup_{t \rightarrow +\infty} \widetilde{\theta}(t, g^{(a)}(t) + x)\right) \] 

\noindent that $\widetilde{\theta}_N^{(a)} = \widetilde{\theta}_N$ takes the same value for each $a \in (0,1)$.

\item Let $\widetilde{u}$ solve Equation~\eqref{eqtildu}, then $\widetilde{\theta}$ solves Equation~\eqref{eqtilth}. This is non-negative, non-decreasing and bounded from above by $N$. Furthermore, $\widetilde{\theta} (t, g^{(a)}(t)) > 0$ (where the inequality is strict), otherwise $\theta (t, x) \equiv 0 \qquad \forall x \leq g^{(a)}(t)$ and hence $\widetilde{u}(t,x) \equiv a \qquad \forall x < g^{(a)}(t)$ and similar contradiction is reached as in the previous part. It follows that, for all $t > 0$,  

\[ \lim_{a \rightarrow 0} \widetilde{\theta}_{x}(t,g^{(a)}(t) ) = \lim_{x \rightarrow +\infty} \widetilde{\theta}_{x}(t,g^{(1/2)}(t) +x) = 0.\]

\noindent It therefore follows that $\lim_{a \rightarrow +\infty} (\log \widetilde{\theta})_x(t,g^{(a)}(t)) = 0$ for all $t > 0$.    

\item This now follows directly from the ergodic theorems proved earlier, taking $a \downarrow 0$.
\end{enumerate}
\qed 

\begin{Th}
\[ \widetilde{\theta}_N = \left\{ \begin{array}{ll} \sqrt{\frac{2\overline{v}}{\kappa}} & N \geq \sqrt{\frac{2\overline{v}}{\kappa}} \\ N & N \leq \sqrt{\frac{2\overline{v}}{\kappa}} \end{array}\right. \]

\noindent and $\overline{g}$ satisfies:

\[ \overline{g} = \left\{ \begin{array}{ll} \sqrt{2\kappa \overline{v}} & N \geq \sqrt{\frac{2\overline{v}}{\kappa}} \\ \frac{\overline{v}}{N} + \frac{\kappa N}{2} & N \leq \sqrt{\frac{2\overline{v}}{\kappa}} \end{array}\right. \]
\end{Th}

\paragraph{Proof} Upper bounds were established by Lemma~\ref{lmmmark}.    Now use:

\begin{equation}\label{eqgtheq} \overline{g} = \frac{\overline{v}}{\widetilde{\theta}_N} + \frac{\kappa \widetilde{\theta}_N}{2}.
\end{equation}

\noindent By considering $\frac{d}{dx} \left(\frac{\overline{v}}{x} + \frac{\kappa x}{2}\right)$ and finding the minimiser, it follows that 

\[ \overline{g} \geq \sqrt{2\kappa \overline{v}}.\]

\noindent Using the upper bound, it therefore follows that, for $N \geq \sqrt{\frac{2\overline{v}}{\kappa}}$, $\overline{g} = \sqrt{2\kappa \overline{v}}$ and hence, from~\eqref{eqgtheq} that for $N \geq \sqrt{\frac{2\overline{v}}{\kappa}}$, $\widetilde{\theta}_N = \sqrt{\frac{2\overline{v}}{\kappa}}$. \vspace{5mm}

\noindent Now consider $N \leq \sqrt{\frac{2\overline{v}}{\kappa}}$. Since $\widetilde{\theta}_N \leq N \leq \sqrt{\frac{2\overline{v}}{\kappa}}$, it follows that $\widetilde{\theta}_N \leq \sqrt{\frac{2 \overline{v}}{\kappa}}$ and hence that

\[ \overline{g} = \frac{\kappa \widetilde{\theta}_N}{2} + \frac{\overline{v}}{\widetilde{\theta}_N} > \kappa \widetilde{\theta}_N. \]

\noindent Let $\widetilde{u}$ satisfy Equation~\eqref{eqtildu} and let $u^*(t,.) = \widetilde{u}(t,g(t) + .)$ From Equation~\eqref{eqtilth}, it follows that $ \theta^*(t,.) := \widetilde{\theta}(t,g(t) + .)$ has a representation

\begin{equation}\label{eqthetrep} \theta^*(t,x) =   N \mathbb{E} \left [ {\bf 1}_{(-\frac{1}{N}\log 2, +\infty)}(Y_{0,t}(x))\exp\left\{ - \int_0^t v(s) u^*(s,Y_{s,t}(x))ds \right\} \right ]\end{equation}

\noindent where

\[ Y_{s,t}(x) = x + (g(t) - g(s)) + (B(t) - B(s)) - \kappa \int_s^t \widetilde{\theta}(r,Y_{r,t}(x))dr.\]

\noindent Let $y = \liminf_{t \rightarrow +\infty} \frac{Y_{0,t}}{t}$, then

\[ y \geq \overline{g} - \kappa \widetilde{\theta}_N > 0\]

\noindent and hence, from~\eqref{eqthetrep} and using the fact that $u^*$ has exponential decay in the space variable,  $\widetilde{\theta}_N = \lim_{x \rightarrow +\infty} \left(\liminf_{t \rightarrow +\infty} \theta^*(t,x)\right) = N$, and $\overline{g} = \frac{\kappa N}{2} + \frac{\overline{v}}{N}$ as required. The theorem is proved. \qed

\section{Noise with Integrable Space Correlation}

Now let $\zeta$ be a Gaussian random field satisfying 
\begin{enumerate} \item $\mathbb{E}\left [ \zeta \right ] \equiv 0$ and 
\item $\mathbb{E} \left [ \zeta(s,x)\zeta(t,y) \right ] = (s \wedge t) \Gamma(x - y)$ where $\Gamma \in C^\infty (\mathbb{R})$ and $\int |\Gamma(z)|dz < +\infty$. 
\end{enumerate} 

\noindent That is, for each $x \in \mathbb{R}$, $\{\zeta(t,x) : t > 0\}$ is a Wiener process with diffusion coefficient $\Gamma(0)$. The processes are correlated, but the correlation decays sufficiently quickly so that the spatial covariance function $\Gamma$ is integrable. Let $u$  satisfy

\begin{equation}\label{eqspace} \left\{ \begin{array}{l} \partial_t u = \left( \frac{\kappa}{2} u_{xx}  + u\left (1 -u \right ) \right) dt + \epsilon u \partial_t \zeta \\ 
             u(0,x) = u_{N,0}:= {\bf 1}_{(-\infty,-\frac{1}{N}\log 2]}(x) + \frac{1}{2}e^{-Nx}{\bf 1}_{(-\frac{1}{N}\log 2, +\infty)}  
           \end{array} \right. \end{equation}
           
\noindent The symbol $\partial_t$ denotes a differential with respect to the `time' variable $t$.  The aim of this section is to study travelling wave front properties of Equation~\eqref{eqspace}; the {\em speed} and the {\em exponential decay}. It turns out that, contrary to situation of the previous section where $\Gamma \equiv 1$,  the asymptotic speed and asymptotic exponential decay of the  travelling front do not depend on $\epsilon$. The proofs of these results given here require integrability of $\Gamma$; also $\Gamma$ must be at least four times differentiable, since heavy use is made of the existence of the second derivatives of the solution in the proofs. The results may hold when $\zeta$ is replaced by `white noise', that is $\Gamma(z) = \delta_0(z)$ (Dirac delta function)  and the stochastic integral taken in the sense of Walsh~\cite{Wa}, but the analysis has not been carried out here and substantially different proofs would be needed. This article gives no information on the white-noise problem. \vspace{5mm}

\noindent The first task is to define a suitable wave marker. Since solutions $u$ to Equation~\eqref{eqspace} do not have the property that $u_x(t,x) \leq 0$ for all $x \in \mathbb{R}$, the definition from Section~\ref{secseason} is difficult to use. A different, deterministic definition of a marker is therefore used for this problem.

\begin{Defn}[Marker]\label{defspcornoisemark} Let $u$ solve Equation~\eqref{eqspace} and let $a^* = \sup_x \left(\liminf_{t \rightarrow +\infty} \mathbb{E} \left [   u(t,x) \right ]  \right)$. Provided $a^* > 0$, then for $a \in (0,a^*)$, the a-marker $g^{(a)} (t)$ is defined as the (deterministic) function such that $ \mathbb{E} \left [   u(t, g^{(a)}(t))\right ] = a$ $\forall t > 0$. If $a^* = 0$ then the marker is not defined.
\end{Defn}

\noindent It is necessary to show that this wave marker is well defined. This will follow from Corollary~\ref{cygamgo}.  \vspace{5mm}

\noindent The speed of the candidate wave marker and the rate of exponential decay beyond the wave marker is given by the following theorem, which shows that the results are exactly the same as those for the deterministic KPP equation (i.e. with $\epsilon = 0$).

\begin{Th}\label{Thcol} Let $u$ solve Equation~\eqref{eqspace} and let $g^{(a)}$ be defined as in Definition~\ref{defspcornoisemark}. Assume that $a^* > 0$. Then for all $a \in (0,a^*)$,

\[ \dot{g}^{(a)} \rightarrow \gamma = \left\{ \begin{array}{ll}   \sqrt{2 \kappa} & N > \sqrt{\frac{2 }{\kappa}} \\ \frac{\kappa N}{2} + \frac{1}{N} & N \leq \sqrt{\frac{2 }{\kappa}} \end{array} \right. \] 

\noindent Furthermore, for all $a \in (0,a^*)$, 

\[ \lim_{x \rightarrow +\infty} \left (\lim_{t \rightarrow +\infty} - \frac{\partial}{\partial x} \mathbb{E} \left [ \log u(t,g^{(a)}(t) + x)\right ] \right ) = \left\{\begin{array}{ll}  \sqrt{\frac{2 }{\kappa}} & N > \sqrt{\frac{2}{\kappa}} \\ N & N \leq \sqrt{\frac{2}{\kappa}}. \end{array}\right. \]
\end{Th}

\noindent Most of the remainder of the section is devoted to establishing results which enable Theorem~\ref{Thcol} to be proved. First, though, Equation~\eqref{eqspace} defines a {\em Feller process}. It is useful to establish the framework for this and indicate the properties that will be used. 

\paragraph{The Solution as a Feller Process} Let $u^{(\gamma)}$ satisfy:

\begin{equation}\label{equg}
\left\{ \begin{array}{l} \partial_t u^{(\gamma)} = \left(\frac{\kappa}{2}u^{(\gamma)}_{xx} + \gamma u^{(\gamma)}_x + u^{(\gamma)} - u^{(\gamma)2}\right) dt + \epsilon u^{(\gamma)}\partial_t \zeta \\ u^{(\gamma)}(0,.) = u_0 \end{array}\right.
\end{equation}

\noindent where $\gamma \in \mathbb{R}$ and $u_0$ is an initial condition, non-negative and uniformly bounded from above. Let 

\begin{equation}\label{eqfnspace} \left\{ \begin{array}{l}  {\cal L}_M = \{ u: \mathbb{R} \rightarrow \mathbb{R}_+ \quad \mbox{Borel measurable}: \sup_x u(x) < M\} \\
{\cal L} = \cup_M {\cal L}_M. \end{array}\right. \end{equation} 

\noindent Let $d$ be the metric over $\mathbb{R}$ defined by:

\[ d(x,y) = \left |\int_x^y \frac{1}{(1+t)^2}dt \right | = \left |\frac{1}{(1+y)} - \frac{1}{(1+x)}\right |. \]

\noindent Let $F : \mathbb{R}_+ \cup \{+\infty\} \rightarrow [0,1]$ be defined by: $F(y) = 1 - e^{-y}$, for $A \in {\cal B}(\mathbb{R})$ let 

\[ A_\alpha = \{x \in \mathbb{R} : \inf_{y \in A} d(x,y) \leq \alpha\}\]

\noindent and let $D$ be the Prokhorov-style metric over ${\cal L}$ defined by:

\begin{eqnarray}\label{eqmetric} \lefteqn{ D(u_1,u_2) = \inf \left \{\alpha: \frac{1}{2}\int_A e^{-|x|}F(u_1(x))dx \leq \alpha + \frac{1}{2}\int_{A_\alpha} e^{-|x|}F(u_2(x))dx \right. }\\&& \left. \qquad \mbox{and} \qquad \frac{1}{2}\int_A e^{-|x|}F(u_2(x))dx \leq \alpha + \frac{1}{2}\int_{A_\alpha} e^{-|x|}F(u_1(x))dx \qquad \forall A \in {\cal B}(\mathbb{R}) \right \}. \nonumber
\end{eqnarray}

\noindent Let $\overline{\cal L}$ denote the completion of ${\cal L}$ under metric $D$. Then $(\overline{\cal L}, D)$ is compact, by the extension of Prokhorov's theorem to finite measures. The space $(\mathbb{R}_+ \cup \{+\infty\},d)$ is clearly compact, using the metric $d$. For a function $u \in \overline{\cal L}$, the measure $\nu_u$ over ${\cal B}(\mathbb{R}_+ \cup \{+\infty\})$ defined by

\[ \nu_u(A) := \frac{1}{2}\int_A e^{-|x|}F(u(x))dx \]

\noindent is clearly of finite variation (the variation is bounded above by $1$). It follows that for any sequence $(u^{(n)})_{n \geq 1}$ of elements in $\overline{{\cal L}}$, there is a convergent subsequence $(u^{(n_j)})_{j \geq 1}$ and limit $u$, in the sense that the measures $(\nu_{n^{(n_j)}})_{j \geq 1}$ converge to a limit $\nu$ such that $\nu(A) = \frac{1}{2}\int_Ae^{-|x|}G(x)dx$. Let $u(x) = \frac{1}{\log(1-G(x)}$, then $\lim_{j \rightarrow +\infty} D(u^{(n_j)},u) = 0$. \vspace{5mm}

\noindent Now let ${\cal P}(\overline{\cal L},D)$ denote the set of probability distributions over $(\overline{\cal L},D)$. Then, since $(\overline{\cal L}, D)$ is compact, it follows that $({\cal P}(\overline{\cal L},D), \rho)$ is compact, by Prokhorov's theorem, where $\rho$ is the Prokhorov metric, defined by:  

\[ \left\{ \begin{array}{l}  \rho(\lambda_1, \lambda_2) = \inf\{\alpha > 0 : \lambda_1 (A) \leq \lambda_2(A_\alpha) + \alpha, \quad \lambda_2(A) < \lambda_1(A_\alpha) + \alpha \qquad \forall A \in {\cal B}(\overline{\cal L})\} \\ A_\alpha = \{ u  \in \overline{{\cal L}} : \inf_{b   \in A} D(u , b ) \leq \alpha\}.\end{array}\right. \]

\noindent Let $\mu^{(T)}$ denote the probability distribution of $u^{(\gamma)}(T,.)$ where $u^{(\gamma)}$ is the solution of~\eqref{eqspace}.  It follows that $(\mu^{(T)})_{T \geq 0}$ has a convergent subsequence and limit $\mu$ which is a probability measure over  $(\overline{\cal L},d)$. The process $u^{(\gamma)}$ defined by~\eqref{eqspace} is Feller. Let $(P_t)_{t \geq 0}$ denote the transition semigroup. Having established existence of a limit $\mu$, it follows that, for a given initial condition $u_0$, this limit is unique and that  $\lim_{t \rightarrow +\infty} \rho (P_t(u_0,.), \mu_N) = 0$. It now follows, in the same way as the previous section, using Theorems 20.17  from Kallenberg~\cite{Kal}  that $(u^{(\gamma)} (t,.))_{t \geq 0}$ is Harris recurrent with supporting measure $\mu$ and hence (Theorem 20.12 page 400 of Kallenberg~\cite{Kal}) it is strongly ergodic.

\begin{Lmm}\label{lmminvmeasuniq} Consider~\eqref{equg} with  non-negative initial conditions. For any fixed $\gamma \in \mathbb{R}$, there is at most one stationary distribution for~\eqref{equg} satisfying 

\[ \inf_x \left(\lim_{t \rightarrow +\infty} \mathbb{E}[u^{(\gamma)}(t,x)]\right) > 0.\]

\noindent This stationary distribution is the same for all $\gamma  \in \mathbb{R}$ and is spatially invariant.  
\end{Lmm}

\paragraph{Proof} Consider two solutions of~\eqref{equg}, with different strictly positive initial conditions. Let $R = \frac{u_1}{u_2}$, then a straightforward application of Itô's formula gives:

\begin{equation}\label{eqRev} \frac{\partial}{\partial t}R = \frac{\kappa}{2}R_{xx} + \left(\gamma +  \frac{\kappa}{2}(\log u_2)_x\right) R_x  + u_2 R(1-R).\end{equation}

\noindent For any local minimum $x_*$ satisfying $R(t,x_*) < 1$, $\left. \frac{\partial}{\partial t}R(t,x) \right |_{x = x_*} > 0$ and for any local maximum $x^*$ satisfying  $R(t,x^*) > 1$, $\left. \frac{\partial}{\partial t}R(t,x)\right |_{x = x^*} < 1$. Suppose $u_1$ and $u_2$ satisfy the hypotheses of the lemma. As $t \rightarrow +\infty$, let $(u_2,R)$ be chosen according to any stationary distribution for the problem, where the evolution of $R$ is defined by~\eqref{eqRev}. Then $R$ will have no local maxima greater than $1$ and no local minima less than $1$. Furthermore, if $\lim_{x \rightarrow +\infty} R(x) = c > 1$, then it follows from~\eqref{eqRev} that $c = +\infty$, contradicting the hypotheses of the lemma; similarly, if $\lim_{x \rightarrow -\infty} R(x) = c > 1$ then $c = +\infty$, again a contradiction. Similarly, $\lim_{x \rightarrow +\infty}R(x) = 0$ or $\lim_{x \rightarrow -\infty}R(x) = 0$ can be excluded by the hypotheses of the lemma. It therefore follows that $R \stackrel{t \rightarrow +\infty}{\longrightarrow} 1$.   

The remaining parts of the lemma follow easily: firstly, with initial condition $u_0 \equiv c > 0$ (equal to a strictly positive constant), it is clear that 1) the invariant measure is spatially invariant and 2) having noted that it is spatially invariant, it is the same for all $\gamma \in \mathbb{R}$, since this is a shift operator. In other words, let $\widetilde{\zeta}^{(\gamma)}(t,x) = \int_0^t \partial_t \zeta(t,x+\gamma t)$. That is, the $\partial_t$ means the differential with respect to the first variable. Then $\widetilde{u}^{(\gamma)}(t,x) = u^{(0)}(t, x+\gamma t)$ where $\widetilde{u}^{(\gamma)}(t,x)$ satisfies: 

\[ \partial_t \widetilde{u}^{(\gamma)} =  \left(\frac{\kappa}{2} \widetilde{u} + \gamma \widetilde{u}^{(\gamma)}_x + \widetilde{u}^{(\gamma)} - \widetilde{u}^{(\gamma)2}\right)dt + \epsilon \widetilde{u}^{(\gamma)} \partial_t \widetilde{\zeta}^{(\gamma)}.\]

\noindent Since $\widetilde{\zeta}^{(\gamma)} \stackrel{(d)}{=} \zeta$, the result follows. \qed \vspace{5mm}

\noindent The following lemma is the most important result for establishing that the `slow-down' from the Itô noise is exactly matched by the `speed-up' from the spatial variance in the noise.

\begin{Lmm}\label{lmmlapfnw}
Let $u^{(\gamma)}$ satisfy~\eqref{equg} with initial condition
\begin{equation}\label{eqic} u_0 = {\bf 1}_{(-\infty,-\frac{1}{N}\log 2]}(x) + e^{-Nx} {\bf 1}_{(-\frac{1}{N}\log 2, +\infty)}(x).
\end{equation}

\noindent Let $v^{(\gamma)} = -\log u^{(\gamma)}$ and $w^{(\gamma)} = v^{(\gamma)}_x$. If $\lim_{t \rightarrow +\infty} \mathbb{E}[u^{(\gamma)}(t,x)] = 0$ for all $x \in \mathbb{R}$, then the limiting log Laplace functional of $w^{(\gamma)}$ satisfies:

\[ {\cal L}(\phi) := \lim_{t \rightarrow +\infty} \log \mathbb{E} \left [e^{\langle \phi, w^ {(\gamma)}(t,.) \rangle}\right ] =  \mu \langle \phi, {\bf 1} \rangle  + \frac{\epsilon^2}{2 \kappa} \int \int \phi(x)\phi(y)\Gamma(x-y)dy dx\]

\noindent where $\langle f,g \rangle := \int_{-\infty}^\infty f(x)g(x)dx$, ${\bf 1}$ denotes the function that is identically equal to $1$ on $\mathbb{R}$ and $\mu \in \mathbb{R}$.
\end{Lmm}

\paragraph{Proof of Lemma~\ref{lmmlapfnw}} Firstly, from the above discussion, that the evolution of  $u$ defines a Feller process and the arguments from which it may be concluded that, for a given initial condition, the distributions of  $u^{(\gamma)}(t,.)$ converge to a unique law, it follows that $\lim_{t \rightarrow +\infty} \mathbb{E}[u(t,x)]$ is well defined. If this limit is equal to $0$ for all $x \in \mathbb{R}$, then the stationary distribution is $u \equiv 0$ with probability $1$. With $v^{(\gamma)} = -\log u^{(\gamma)}$ and $w^{(\gamma)} = v_x^{(\gamma)}$, Itô's formula gives:

\begin{equation}\label{equv}   
\partial_t w^{(\gamma)} = \left(\frac{\kappa}{2} w_{xx}^ {(\gamma)} - \frac{\kappa}{2}(w^{(\gamma)2})_x  + \gamma w^{(\gamma)}_x - u^{(\gamma)} w^{(\gamma)}\right)dt - \epsilon \partial_t \zeta_x
 \end{equation}
           
\noindent Here $\zeta_x$ denotes the derivative of $\zeta$ with respect to the $x$ variable. Under the assumption that $\mathbb{E}[u(t,x)] \rightarrow 0$ for each $x \in \mathbb{R}$, it follows that, asymptotically, $w^{(\gamma)}$  satisfies:

\begin{equation}\label{equvasym} 
\partial_t w^{(\gamma)} = \left(\frac{\kappa}{2} w_{xx}^{(\gamma)} - \frac{\kappa}{2}(w^{(\gamma)2})_x + \gamma w_x^{(\gamma)} \right)dt - \epsilon \partial_t \zeta_x.
\end{equation}

\noindent Furthermore, the  evolution of $w^{(\gamma)}$ defines  a strongly ergodic Feller process whose stationary distribution satisfies the second equation of~\eqref{equvasym}. By ergodicity, the stationary distribution will be spatially invariant. 

Let ${\cal H}(t,\phi) = \mathbb{E} \left [ e^{\langle \phi,    w(t)\rangle}\right ]$, where $\phi \in {\cal S}$ defined by~\eqref{eqtfsp}:

\begin{equation}\label{eqtfsp} {\cal S} = \left \{ \phi : \int_{-\infty}^\infty \left(|\phi_x(x)| + |\phi(x)|\right ) dx < +\infty, \qquad \exists x^* : \phi(x^*+y) = \phi(x^* -y) \quad \forall y \in \mathbb{R} \right \} \end{equation}

\noindent and $\langle \phi, \psi \rangle = \int_{-\infty}^\infty \phi(x) \psi (x) dx$.  This is the space of functions over which it is straightforward to compute the limit of the log Laplace functional as $t \rightarrow +\infty$. This will be sufficient to characterise  the limiting distribution and therefore the result is true for all $\phi : \int_{-\infty}^\infty |\phi_x(x)| + |\phi(x)| dx < +\infty$. 

An application of Itô's formula gives:

\begin{eqnarray*} \frac{\partial}{\partial t} {\cal H}(t,\phi) &=& \frac{\kappa}{2} \left \langle \phi_{xx}, \mathbb{E}\left [e^{\langle \phi, w(t) \rangle}w(t)\right ] \right \rangle - \gamma \left \langle \phi_x, \mathbb{E} \left [ e^{\langle \phi, w(t) \rangle}w(t)\right ] \right \rangle - \left \langle \phi, \mathbb{E} \left [e^{\langle \phi, w(t)\rangle}u(t)  w(t) \right ] \right \rangle \\&& + \frac{\kappa}{2} \left \langle \phi_x, \mathbb{E} \left [e^{\langle \phi,w(t)\rangle} w(t)^2\right ] \right \rangle - \frac{\epsilon^2}{2} {\cal H}(t,\phi)\int \int \phi(x)\phi(y) \Gamma^{\prime \prime}(x-y)dy dx 
\end{eqnarray*}

\noindent where $\Gamma^{\prime \prime}$ denotes the second derivative of $\Gamma$. The functional equation may be written

\begin{eqnarray*}  \frac{\partial}{\partial t}{\cal H}(t,\phi) &=& \frac{\kappa}{2} \int \phi_{xx}(x) \frac{\partial}{\partial \delta (x)}{\cal H}(t,\phi) dx - \gamma \int \phi_x(x)\frac{\partial}{\partial \delta(x)}{\cal H}(t,\phi) dx  - \left \langle \phi, \mathbb{E} \left [e^{\langle \phi, w(t)\rangle}u(t)  w(t) \right ] \right \rangle \\&& + \frac{\kappa}{2}\int \phi_x(x)\frac{\partial^2}{\partial \delta(x)^2} {\cal H}(t,\phi) dx - \frac{\epsilon^2}{2}{\cal H}(t,\phi)\int \int \phi(x)\phi(y)\Gamma^{\prime \prime}(x-y)dy dx.
\end{eqnarray*}

\noindent where $\frac{\partial}{\partial \delta(x)}$ denotes the {\em functional} derivative with respect to $\delta(x)$, which is defined as: 

\[ \frac{\partial}{\partial \delta(x)}{\cal H}(t,\phi) = \lim_{\epsilon \rightarrow 0}\frac{{\cal H}(t, \phi + \epsilon \delta (x)) - {\cal H}(t,\phi)}{\epsilon};\]

\noindent $\delta(x)$ is the dirac delta function with unit mass at point $x \in \mathbb{R}$. Let ${\cal H}(\phi) = \lim_{t \rightarrow +\infty} {\cal H}(t,\phi)$. It follows that ${\cal H}$ satisfies:

\begin{equation}\label{eqhphi}\begin{array}{r} \int \phi(x)\left\{ \frac{\kappa}{2} \frac{\partial^2}{\partial x^2} \left(\frac{\partial {\cal H} (\phi)}{\partial \delta(x)}\right)  - \frac{\kappa}{2} \frac{\partial}{\partial x} \left(\frac{\partial^2 {\cal H}(\phi) }{\partial \delta(x)^2} \right) + \gamma \frac{\partial}{\partial x}\left(\frac{\partial {\cal H}(\phi)}{\partial \delta(x)}\right)\right\} dx \\ = \frac{\epsilon^2}{2}{\cal H}(t,\phi) \int \int \phi(x)\phi(y)\Gamma^{\prime \prime}(x-y)  dydx. \end{array} 
\end{equation}

\noindent Let ${\cal L} = \log {\cal H}$, then it is straightforward from~\eqref{eqhphi}, by firstly integration by parts so that the spatial derivatives are taken with respect to $\phi$, that:

\begin{equation}\label{eqloglap} \begin{array}{l} 
\int \phi(x) \left\{\frac{\kappa}{2}\frac{\partial^2}{\partial x^2} \left(\frac{\partial {\cal L}(\phi)}{\partial \delta(x)}\right) - \frac{\kappa}{2}\frac{\partial}{\partial x} \left(\frac{\partial {\cal L}(\phi)}{\partial \delta(x)}\right)^2 - \frac{\kappa}{2}\frac{\partial}{\partial x} \left(\frac{\partial^2 {\cal L}(\phi)}{\partial \delta(x)^2}\right) + \gamma \frac{\partial}{\partial x}\left(\frac{\partial  {\cal L}(\phi)}{\partial \delta(x)}\right)  \right\}dx \\ \hspace{50mm}  = \frac{\epsilon^2}{2} \int \int \phi(x) \phi(y) \Gamma^{\prime \prime}(x-y) dy   dx.  \end{array} 
\end{equation}

 \noindent Set 

\[{\cal L}_0(\phi) :=   \frac{\epsilon^2}{2\kappa}  \int \int \phi(x)\phi(y)\Gamma(x-y)dydx \]

\noindent then ${\cal L}_0$ is a solution to~\eqref{eqloglap}. This is seen as follows: direct computation gives:

\[ \frac{\partial {\cal L}_0(\phi)}{\partial \delta(x)} =  \frac{\epsilon^2}{\kappa} \int \phi(y)\Gamma(x-y)dy\]

\noindent and 

\[ \frac{\partial^2 {\cal L}_0(\phi)}{\partial \delta(x)^2} = \frac{\epsilon^2}{\kappa}\Gamma(0).\]

\noindent It follows that

\[ \frac{\kappa}{2} \int \phi(x)\frac{\partial^2}{\partial x^2} \left(\frac{\partial {\cal L}_0 (\phi) }{\partial \delta(x)}\right)dx = \frac{\epsilon^2}{2}\int \int \phi(x)\phi(y) \Gamma^{\prime \prime}(x-y)dx dy  \]

\noindent  and

\[ \int \phi(x)\frac{\partial}{\partial x}\left(\frac{\partial^2 {\cal L}_0(\phi)}{\partial \delta(x)^2}\right)dx = \int \phi(x)\frac{\partial}{\partial x}\left(\frac{\partial  {\cal L}_0(\phi)}{\partial \delta(x) }\right)dx = 0.\]

\noindent For the remaining terms: 

\[
  \int\phi(x) \frac{\partial}{\partial x} \left(\frac{\partial {\cal L}_0(\phi)}{\partial \delta(x)} \right)^2 dx =   \frac{\epsilon^2}{\kappa^2} \int \int \int \phi(x)\phi(y)\phi(z) \frac{\partial}{\partial x} (\Gamma(x-y)\Gamma(x-z))dydzdx
\]

\noindent Since $\Gamma(x-y) = \Gamma(y-x)$, it follows that: 

\[ \int \int \phi(x)\phi(y)\Gamma^\prime(x-y)dy dx = \int \int \phi(x) \phi^\prime(y)\Gamma(x-y)dy dx = -\int \int \phi^\prime(x) \phi(y)\Gamma(x-y) dy dx = 0\]

\noindent Now set

\[ I = \int \int \int \phi^\prime (x)\phi(y)\phi(z) \Gamma(x-y)\Gamma(x-z)  dy dz dx\]

\noindent Since $\phi$ is {\em symmetric} around a point $x^*$ (i.e. $\phi(x^* + x) = \phi(x^* -x)$ for all $x \in \mathbb{R}$), it follows directly, using

\[ I = \int \int \int \phi^\prime (x^* + x) \phi(x^*+y)\phi(x^*+z)\Gamma(x-y) \Gamma(x-z) dy dz dx\]

\noindent together with $\phi^\prime(x^*+x) = -\phi^\prime(x^*-x)$ that 

\[ I = -I = 0.\]

\noindent Therefore ${\cal L}_0$ is a solution.
\vspace{5mm}

\noindent Now consider solutions of the form: ${\cal L}$ defined such that ${\cal L}(\phi) = {\cal L}_0(\phi) + {\cal E}(\phi)$. Then  ${\cal E}(\phi)$ satisfies:

\begin{equation}\label{eqesat} \left. \begin{array}{r} \int \phi(x) \left\{ \frac{\kappa}{2} \frac{\partial^2}{\partial x^2} \left(\frac{\partial {\cal E} (\phi)}{\partial \delta(x)}\right) - \frac{\kappa}{2} \frac{\partial}{\partial x} \left(\left(2\frac{\partial {\cal L}_0 (\phi)}{\partial \delta(x)} + \frac{\partial {\cal E}(\phi)}{\partial \delta(x)}\right)\left(\frac{\partial {\cal E}(\phi)}{\partial \delta(x)}\right)\right) \right. \\   -  \left. \frac{\kappa}{2} \frac{\partial}{\partial x} \left(\frac{\partial^2 {\cal E}(\phi)}{\partial \delta(x)^2} \right) + \gamma \frac{\partial}{\partial x}\left(\frac{\partial {\cal E}(\phi)}{\partial \delta(x)}\right) \right\}dx = 0. \end{array}\right.
\end{equation}

\noindent Let $\Gamma \phi (x) := \int \Gamma(x-y)\phi(y)dy$, then~\eqref{eqesat} may be rewritten as:

\begin{equation}\label{eqsat2} \begin{array}{r}
\int \phi(x) \left\{ \frac{\kappa}{2} \frac{\partial^2}{\partial x^2} \left(\frac{\partial {\cal E}(\phi)}{\partial \delta(x)}\right) - \epsilon^2 \frac{\partial}{\partial x} \left((\Gamma \phi)(x) \frac{\partial {\cal E}(\phi)}{\partial \delta(x)}\right) - \frac{\kappa}{2} \frac{\partial}{\partial x} \left(   \frac{\partial {\cal E}(\phi)}{\partial \delta(x)}\right)^2 \right. \\ \left. - \frac{\kappa}{2} \frac{\partial}{\partial x} \left(\frac{\partial^2 {\cal E}(\phi)}{\partial \delta(x)^2} \right) + \gamma \frac{\partial}{\partial x}\left(\frac{\partial {\cal E}(\phi)}{\partial \delta(x)}\right)\right\} = 0. \end{array}
\end{equation}

\noindent Since $\phi \in {\cal S}$ (and is therefore symmetric around some point $x^*$),   it follows that  $\frac{\partial {\cal E}(\phi)}{\partial \delta(x)}$ and $\frac{\partial^2 {\cal E}(\phi)}{\partial \delta(x)^2}$ are symmetric around the same point and therefore ~\eqref{eqsat2} reduces to:

\[ \int \phi(x) \frac{\partial^2}{\partial x^2} \left(\frac{\partial {\cal E}(\phi)}{\partial \delta(x)}\right) dx = 0.\]

\noindent Now, a general expression for ${\cal E}(\phi)$ is:

\[ {\cal E}(\phi) = E_0 + \sum_{n=1}^\infty \int\ldots\int E_n(x_1,\ldots,x_n)\phi(x_1)\ldots\phi(x_n)dx_1\ldots dx_n,\]

\noindent from which (using exchangeability of the variables in $E_n(x_1, \ldots, x_n)$):

\[ \int \phi(x)\frac{\partial^2}{\partial x^2} \left(\frac{\partial {\cal E}(\phi)}{\partial \delta (x)}\right)dx = \sum_{n=1}^\infty \frac{1}{n}\int\ldots\int \phi(x_1)\ldots \phi(x_n)  \Delta E_n(x_1,\ldots, x_n) dx_1 \ldots x_n.\]

\noindent where $\Delta = \frac{\partial^2}{\partial x_1^2} + \ldots + \frac{\partial^2}{\partial x_n^2}$ is the Laplacian in $\mathbb{R}^n$. Let 

\[ C_n (\phi) = \frac{1}{n}\int \ldots \int E_n(x_1,\ldots,x_n)\phi(x_1)\ldots\phi(x_n)dx_1\ldots dx_n.\]

\noindent  Since this holds for all $\phi \in {\cal S}$, it follows that for all $\alpha \in \mathbb{R}$, $\sum_n \alpha^n C_n(\phi) = 0$, from which it follows that $C_n(\phi) = 0$ for each each  $n$, $\Delta E_n \equiv 0$. Using $E_n(x_1, \ldots, x_n) = E_n(z+x_1, \ldots, z+x_n)$ for all $z \in \mathbb{R}$, $(x_1, \ldots, x_n) \in \mathbb{R}^n$, it follows that $E_n (x_1, \ldots, x_n) = E_n$, a constant.

 From this,  it follows that ${\cal E}(\phi) \equiv F(\langle \phi, {\bf 1} \rangle)$ for some function $F$, where $\langle \phi, {\bf 1} \rangle := \int_{-\infty}^\infty \phi(x)dx$. Hence it follows that any solution satisfies:

\[ {\cal L}(\phi) = F\left(\langle \phi, {\bf 1} \rangle\right)  +  \frac{\epsilon^2}{2\kappa} \int \int \phi(x) \phi(y) \Gamma(x-y) dy dx \qquad F : \mathbb{R} \rightarrow \mathbb{R}.\]
 
 \noindent Since ${\cal L}$ is a log-Laplace functional, it is convex; for all $\lambda \in (0,1)$ and all $\phi,\psi \in {\cal S}$, ${\cal L}(\lambda \phi + (1-\lambda)\psi) \leq \lambda{\cal L}(\phi) + (1-\lambda){\cal L}(\psi)$. If $\phi,\psi \in {\cal S}$, then $\phi_n, \psi_n$ defined by $\phi_n(x) = \frac{1}{n} \phi\left (\frac{x}{n} \right )$ and $\psi_n(x) = \frac{1}{n}\psi_n \left (\frac{x}{n} \right ) \in {\cal S}$. Since $F \left(\langle \phi_n, {\bf 1} \right ) = F\left(\langle \phi, {\bf 1}\rangle \right )$ for all $n \geq 1$ and (clearly) $\lim_{n \rightarrow +\infty} \int \int \phi_n(x)\phi_n(y)\Gamma(x-y)dydx \stackrel{n \rightarrow +\infty}{\longrightarrow} 0$ (similarly for $\psi_n$ and $\lambda\phi_n + (1-\lambda)\psi_n$), it follows that $F$ is convex.

From the form of the Laplace functional, it follows that for $w$ distributed according to the stationary distribution, $w(x) = \widetilde{w} + \widehat{w}$ where $\widetilde{w}$ has log Laplace functional $F(\langle \widetilde{w}, \phi \rangle)$. It follows that $\widetilde{w}$ is a random variable (i.e. does not depend on $x$). 

Now, from strong ergodicity of $w$ (which may be established using the same machinery as strong ergodicity of $u$), it follows, since $\widetilde{w}$ is a random variable (which does not depend on $x$) that $F(\langle \phi, {\bf 1} \rangle) ) = \mu \langle \phi, {\bf 1} \rangle$ for some constant $\mu \in \mathbb{R}$. 
The lemma has now been proved. \qed 
 
 \begin{Lmm} \label{lmmlimvar} Let $u$ satisfy~\eqref{equg} with initial condition~\eqref{eqic}. Suppose that the stationary distribution $\mu_{N,\gamma}$ satisfies $\lim_{x \rightarrow +\infty}\mathbb{E}_{\mu_{N,\gamma}}[u(x)] = 0$, where $u$ has distribution  $\mu_{N,\gamma}$. Let $w = -\frac{\partial}{\partial x} \log u(t,x)$ and, for $\phi \in {\cal S}$, let ${\cal L}(\phi) = \lim_{z \rightarrow +\infty} \left (\lim_{t \rightarrow +\infty} \log \mathbb{E} \left [ \exp\left\{ \langle w(t),\phi(.+z) \right\}\right ] \right)$ then
 
 \[{\cal L}(\phi) = \mu \langle {\bf 1}, \phi \rangle + \frac{\epsilon^2}{2\kappa} \int \int \phi(x)\phi(y)\Gamma(x-y)dydx \]
 
 \noindent where $\mu \in \mathbb{R}$. 
 \end{Lmm}
 
 \paragraph{Proof} The proof is the same as that for Lemma~\ref{lmmlapfnw}; the only issue is to establish that  $F(\langle \phi, {\bf 1} \rangle) = \mu \langle \phi, {\bf 1}\rangle$ for some  $\mu \in \mathbb{R}$. In Lemma~\ref{lmmlapfnw}, this followed directly from the fact that the evolution defined a Feller process which, following standard results, was strongly ergodic. 
 
  The same arguments apply here; $\lim_{z \rightarrow +\infty} {\cal L}(\phi(.+z))$ defines the invariant measure of a Feller process with transition semigroup defined by ${\cal P}_t^*(w,A) = \lim_{a \rightarrow +\infty} {\cal P}_t(w(.+a),A\circ \theta_a)$ $\theta_a$ denotes the shift operator, $\theta_a f(.) = f(a+.)$ for functions $f$, where ${\cal P}_t$ is the transition semigroup of the pair $(u,w)$. The same arguments as Lemma~\ref{lmmlapfnw} now give the result. \qed\vspace{5mm}
 
 \noindent The following result is a corollary of Lemma~\ref{lmmlapfnw}.

 \begin{Cy}\label{cygamgo} Let $u$ satisfy~\eqref{equg} with $\gamma = 0$ and initial condition $u_0$ given by~\ref{eqic}, with  $N \geq 0$. Then $\lim_{t \rightarrow +\infty} \mathbb{E}[u(t,x)] > 0$ for each $x \in \mathbb{R}$. 
 \end{Cy}

 \paragraph{Proof} Let $v = -\log u$, then $v$ satisfies:
 
 \[ \partial_t v = \left(\frac{\kappa}{2} v_{xx} - \frac{\kappa}{2}w^2 - 1 + u + \frac{\epsilon^2}{2}\Gamma(0)\right) dt - \epsilon \partial_t \zeta. \]
 
\noindent Lemma~\ref{lmmlapfnw} gives directly that:
 
\[ \lim_{t \rightarrow +\infty} \mathbb{E} \left [ w(t,x)w(t,y) \right ] = \mu^2 + \frac{\epsilon^2}{\kappa} \Gamma(x-y).\]

\noindent Let $V(t,x) = \mathbb{E} \left [ v(t,x) \right ]$. Then, under the assumption that $\lim_{t \rightarrow +\infty} \mathbb{E}[u(t,x)] = 0$ for each $x \in \mathbb{R}$, it follows that: 

\[ \lim_{t \rightarrow +\infty} \frac{\partial}{\partial t} V = -\frac{\kappa}{2} \mu^2 -1.\]

\noindent  It follows that $- V(t,x) = \mathbb{E}[\log u(t,x)] \stackrel{t \rightarrow +\infty}{\longrightarrow} + \infty$ and since, by Jensen,  $\mathbb{E}[\log u(t,x)] \leq \log \mathbb{E}[u(t,x)]$, this   contradicts the fact that $\mathbb{E}[u(t,x)]\stackrel{t \rightarrow +\infty}{\longrightarrow} 0$.  Corollary~\ref{cygamgo} is proved. \qed

\begin{Cy}\label{cyuinvmeas}  Let $u$ solve~\eqref{eqspace}. Let $a^*$ be defined in Definition~\ref{defsmark}. Then,   $u$ has a unique invariant measure which is independently  of the value of $N$ in the initial condition and satisfies  $\lim_{t \rightarrow +\infty} \mathbb{E}[u(t,x)] = a^*$ for each $x \in \mathbb{R}$. Furthermore,  $a^* > 0$.  
\end{Cy}

\paragraph{Proof} This follows from Corollary~\ref{cygamgo} together with Lemma~\ref{lmminvmeasuniq}. \qed

\begin{Lmm}\label{lmmgotlim}
 For each $a \in (0,a^*)$,

\begin{equation}\label{eqbdaprior} 0 \leq \liminf_{t \rightarrow +\infty}\frac{g^{(a)}}{t} \leq \limsup_{t \rightarrow +\infty} \frac{g^{(a)}(t)}{t} \leq \left\{ \begin{array}{ll} \sqrt{2\kappa} & N > \sqrt{\frac{2}{\kappa}} \\ \frac{\kappa N}{2} + \frac{1}{N} & N \leq \sqrt{\frac{2}{\kappa}}. \end{array}\right. \end{equation}

\end{Lmm}

\paragraph{Proof of Lemma~\ref{lmmgotlim}} The lower bound is a direct consequence of Corollary~\ref{cyuinvmeas}.  For the upper bound,  $U^{(a)}$ defined by $U^{(a)}(t,x) = \mathbb{E}[u(t,x+g^{(a)}(t))]$ satisfies:

\[ U_t^{(a)} = \frac{\kappa}{2} U^{(a)}_{xx} + \dot{g}^{(a)}U^{(a)}_x + U^{(a)} - \mathbb{E}[u^{2}](t,g^{(a)}(t) + .).\]

\noindent Exactly the same arguments as those for Lemma~\ref{lmmmark} give the a-priori upper bounds on $g^{(a)}$, since the upper bound is obtained by ignoring the nonlinear term, in this case $\mathbb{E}[u^2](t,g^{(a)}(t) + x)$. Let $\gamma^{(a)} = \limsup_{t \rightarrow +\infty} \frac{g^{(a)}(t)}{t}$;~\eqref{eqbdaprior} follows and Lemma~\ref{lmmgotlim} is proved. \qed

\begin{Lmm}\label{lmmlogabd} Let $u$ satisfy~\eqref{eqspace} with $g^{(a)}$ and $a^*$ from Definition~\ref{defspcornoisemark}. Then, for any  $a \in (0,a^*)$,  
  \[ \limsup_{t \rightarrow +\infty} \frac{g^{(a)}(t)}{t} = \liminf_{t \rightarrow +\infty} \frac{g^{(a)}(t)}{t} = \gamma^{(a)}\]
 for a constant $\gamma^{(a)} \geq 0$.
\end{Lmm}

\paragraph{Proof} For fixed $N > 0$, consider the family of equations indexed by $\gamma$: 

\begin{equation}\label{eqgamfam} \left\{ \begin{array}{l} \partial_t u^{(\gamma)}  = \left ( \frac{\kappa}{2}u_{xx}^{(\gamma)} + \gamma u_x^{(\gamma)} + u^{(\gamma)} - u^{(\gamma)2} \right )dt  + \epsilon u^{(\gamma)} \partial \zeta \\ u(0,.) = {\bf 1}_{(-\infty,-\frac{1}{N}\log \frac{1}{a}]}(x) + \frac{1}{2}e^{-Nx}{\bf 1}_{( \frac{1}{N}\log \frac{1}{a} , +\infty)}(x) \end{array}\right. \end{equation}

\noindent This defines a Feller process and, by similar arguments as before, the laws of $u^{(\gamma)}(t,.)$ converge to a stationary distribution $\mu_{N,\gamma}$. Let $e_\gamma = \mathbb{E}_{\mu_{N,\gamma}}[u(0)]$; the expected value of $u(0)$ under the law $\mu_{N,\gamma}$. Then, since $u^{(\gamma)}(t,.) = u^{(0)}(t,\gamma t + .)$, it follows that $e_\gamma$ is decreasing in $\gamma$. A gentle modification of the proof of Lemma~\ref{lmmgotlim} gives that, for any $a  \in (0,a^*)$ (open interval),  $e_\gamma = 0$ for $\gamma > \limsup_{t \rightarrow +\infty}\frac{g^{(a)}(t)}{t}$. Furthermore,  $e_0 = a^*$. 

There is a unique $\gamma : e_{\gamma -} \geq a \geq e_{\gamma +}$. Denote this by $\gamma^{(a)}$. Otherwise, if there is an interval $(\gamma_-,\gamma_+)$ such that $e_\beta = a$ for all $\beta \in (\gamma_-,\gamma_+)$, this implies that for such $\beta$, Equation~\eqref{eqgamfam}, with $\gamma = \beta$, has a spatially invariant stationary distribution satisfying $\mathbb{E}[u(0)] = a < a^*$ which is a contradiction. It therefore follows, almost directly from the definition, that $\lim_{t \rightarrow +\infty} \frac{g^{(a)}(t)}{t} = \gamma^{(a)}$, since clearly for any $\epsilon > 0$, $\gamma^{(a)}-\epsilon < \liminf_{t \rightarrow +\infty} \frac{g^{(a)}(t)}{t} \leq \limsup_{t \rightarrow +\infty} \frac{g^{(a)}(t)}{t} < \gamma^{(a)} + \epsilon$.  \qed

\begin{Lmm}\label{lmmgameq} Let $\gamma^{(a)} =  \lim_{t \rightarrow +\infty} \frac{g^{(a)}}{t}$. Then, for any given $N$, $\gamma^{(a)} = \gamma^{(b)} = \gamma$ for all $0 < a < b < a^*$.
\end{Lmm}

\paragraph{Proof} Assume not, then $\gamma^{(a)} > \gamma^{(b)}$. Let $c \in (a,b)$ and consider~\eqref{eqgamfam} with $\gamma = \gamma^{(c)}$. Then $\lim_{x \rightarrow +\infty} \left(\lim_{t \rightarrow +\infty}\mathbb{E}[u^{(\gamma^{(c)})}(t,x)]\right) \geq a$, hence $u^{(\gamma^{(c)})}$ has stationary distribution with expected value $a^*$, contradicting the fact that its expectation is less than $b < a^*$. \qed \vspace{5mm}

\noindent From now on, let $\gamma = \gamma^{(a)}$ for all $a  \in (0,a^*)$. Note that this value depends on $N$; when it is necessary to make the dependence explicit, $\gamma_N$ will be used.

\begin{Lmm}\label{lmmlimgamlimmeas} Consider Equaation~\eqref{equg} with initial condition~\eqref{eqic} (which depends on $N$). 
For each $a \in (0,a^*)$, $\lim_{t \rightarrow +\infty}\dot{g}^{(a)} = \gamma$, in the sense that the distributions of $u(t,g^{(a)}(t) +  .)$ converge to a stationary distribution $\mu_{N,a}$ such that  $\mathbb{E}_{\mu_{N,a}} [u(0)] = a$ and $\mu_{N,a}$ is an invariant measure for the evolution defined by:

\[ \partial_t u = \left( \frac{\kappa}{2}u_{xx} + \gamma  u_x + u - u^2\right) dt + \epsilon u \partial_t \zeta. \]
\end{Lmm}

\paragraph{Proof} Let $\mu_{N,a}(t)$ denote the laws of $u(t,g^{(a)}(t) + .)$. By compactness, there is a convergent subsequence $(\mu_{N,a}(t_n))_{n \geq 1}$ and a limit point $\mu_{N,a}$. Consider a subsequence such that $\lim_{n \rightarrow +\infty}(t_{n+1} - t_n) = +\infty$. Such a subsequence may be extracted without loss of generality. Let $\gamma (n) = \frac{g^{(a)}(t_{n+1}) - g^{(a)}(t_n)}{t_{n+1} - t_n}$. Then the evolution on the time interval $(t_n,t_{n+1}]$ given by:

\[ \partial_t u = \left (\frac{\kappa}{2}u_{xx} + \gamma (n) u_x + u - u^2\right)dt + \epsilon u \partial_t \zeta \]

\noindent gives precisely $u(t_n,g^{(a)}(t_n) + .)$ at the time point $t_n$ for each $n$. Since (by contruction) $\gamma (n) \rightarrow  \gamma$, it follows that, in the limit, the evolution is given by:

\[ \partial_t u = \left (\frac{\kappa}{2}u_{xx} + \gamma  u_x + u - u^2\right)dt + \epsilon u \partial_t \zeta \]

\noindent and that the measure $\mu_{N,a}$ is an invariant measure for this evolution, with expected value $\mathbb{E}_{\mu_{N,a}}[u(0)] = a$. \qed

\begin{Lmm}
 \begin{itemize}
\item \begin{equation}\label{eqvbound} \log \frac{1}{a} \leq -  \mathbb{E} \left [ \log u(t,g^{(a)}(t))\right ] \qquad \forall a \in (0,a^*) \end{equation}

\item 
\begin{equation}\label{eqvlb} 
\lim_{M \rightarrow +\infty} \sup_{t \in \mathbb{R}_+} \mathbb{P}\left(-\log u(t,g^{(a)}(t)) > M\right) = 0 \qquad \forall a \in (0,a^*)
\end{equation}
\end{itemize}
\end{Lmm}

\paragraph{Proof} For the first part, Jensen's inequality gives:

\[ -\mathbb{E}[\log u(t,g^{(a)}(t))] \geq -\log \mathbb{E}[u(t,g^{(a)}(t))] = \log \frac{1}{a}.\]

\noindent For the second part, let $\mu_{a,N}$ denote the limiting distribution of $u(t, g^{(a)}(t) + .)$. Let $u_0$ be a function chosen according to the distribution $\mu_{a,N}$ and let $\gamma = \lim_{t \rightarrow +\infty} \dot{g}^{(a)}(t)$. Then, for each $t > 0$, $u^{(\gamma)} (t,.)$ also has distribution $\mu_{a,N}$, where $\gamma = \lim_{t \rightarrow +\infty} \dot{g}(t)$.

\noindent Note that the solution $u^{(\gamma)}$ may be written:

\[ u^{(\gamma)}(t,x) = \mathbb{E} \left [u_0(X_{0,t})\exp\left\{\int_0^t (1 - u^{(\gamma)}(s,X_{s,t}))ds + \epsilon \int_0^t \partial_s \zeta(s,X_{s,t}) - \frac{\epsilon^2}{2}\Gamma(0) t\right\} \right ]. \]

\noindent Here 

\[X_{s,t}(x) = x + (B_t - B_s) + \gamma (t-s) \]

\noindent where $B$ is a Wiener process, independent of $\zeta$ and $u_0$, with diffusion coefficient $\kappa$ and the path integral is standard: $\int_0^t \partial_s \zeta(s,X_{s,t}) = \lim_{\mbox{mesh} \rightarrow 0} \sum (\zeta_{s_{i+1}} - \zeta_{s_i})(X_{s_i,t})$. This is a basic Kacs representation. 

\noindent Now, $\int_0^t \partial_s \zeta_s(X_{s,t})$ is a centred Gaussian random variable with variance $\frac{\Gamma(0)}{2}t$, while $\mathbb{E}[u^{(\gamma)}(s,X_{s,t})]  \leq a^* \leq 1$ and $\mathbb{E}[u^{(\gamma )2}(s,X_{s,t})] < a^* \leq 1$. It therefore follows that if $u^{(\gamma)}(t,0) = 0$, then $u_0 \equiv 0$, contradicting the fact that $\mathbb{E}[u_0(0)] = a$ (where, of course, $\zeta$ is independent of $u_0$).

This is established as follows: $u_0(0)$ is the limit of random variables $u(t,g^{(a)}(t))$ satisfying 

\[ \mathbb{E}[u(t,g^{(a)}(t))] = a \qquad \mathbb{E}[u^2(t,g^{(a)}(t))] < a. \] 

\noindent The uniform bound in $L^2$ gives that the sequence is uniformly integrable. It follows that if $u_0 \equiv 0$, then $u(t,g^{(a)}(t))$ converges in $L^1$ norm to $0$, contradicting the fact that $\mathbb{E}[u(t,g^{(a)}(t))] = a$ for all $t > 0$. The result follows.  
\qed

\paragraph{Proof of Theorem~\ref{Thcol}} Let $\tilde{u}^{(a)}(t,.) = u(t,g^{(a)}(t) + .)$. Let $\widetilde{v} = - \log \widetilde{u}$ and $\widetilde{w} = \widetilde{v}_x$.  By Itô's formula,

\begin{equation}\label{equvw} \left\{ \begin{array}{l} \partial_t \widetilde{u}^{(a)} = \left(\frac{\kappa}{2} \widetilde{u}^{(a)}_{xx} + \dot{g}^{(a)}(t)\widetilde{u}^{(a)}_x + \widetilde{u}^{(a)} - \widetilde{u}^{(a)2} \right)dt + \epsilon \widetilde{u}^{(a)} \partial_t \widetilde{\zeta} \\ \partial_t \widetilde{v}^{(a)} = \left(\frac{\kappa}{2}\widetilde{w}^{(a)}_x - \frac{\kappa}{2}\widetilde{w}^{(a)2} + \dot{g}^{(a)}\widetilde{w}^{(a)} - 1 + \widetilde{u}^{(a)} + \frac{\epsilon^2}{2}\Gamma(0)\right)dt - \epsilon \partial_t \widetilde{\zeta} \\ 
\partial_t \widetilde{w}^{(a)}_t =  \left(\frac{\kappa}{2}\widetilde{w}^{(a)}_{xx} - \frac{\kappa}{2}(\widetilde{w}^{(a)})^2_x + \dot{g}^{(a)} \widetilde{w}_x - \widetilde{u}^{(a)} \widetilde{w}^{(a)} \right)dt - \epsilon \partial_t \widetilde{\zeta}_x
           \end{array}\right. \end{equation}

\noindent Let $\gamma = \gamma^{(a)}$ of Lemma~\ref{lmmlogabd} and recall that $\gamma^{(a)}$ are equal for all $a \in (0,a^*)$ by Lemma~\ref{lmmgameq} (they do not depend on $a$, but they do depend on $N$). From Lemma~\ref{lmmlimgamlimmeas} and Equation~\eqref{eqvlb}, it follows from the equation for $\widetilde{v}^{(a)}$ that $\lim_{t \rightarrow +\infty} \frac{v^{(a)}(t,x)}{t} = 0$ in probability and hence, taking $\lim_{t \rightarrow +\infty} \frac{1}{t}\int_0^t$ for each term in the second of~\eqref{equvw} together with strong ergodicity gives: 

\[0 = \frac{\kappa}{2}\mathbb{E}_{\mu_{N,a}}[w(x)]_x - \frac{\kappa}{2} \mathbb{E}_{\mu_{N,a}}[w^2(x)] + \gamma \mathbb{E}_{\mu_{N,a}}[w(x)] - 1 + \mathbb{E}_{\mu_{N,a}}[u(x)] + \frac{\epsilon^2}{2}\Gamma(0)\]

\noindent where $\mu_{N,a}$ is the stationary distibution defined in Lemma~\ref{lmmlimgamlimmeas}, $u \sim \mu_{N,a}$ and $w = -\frac{\partial}{\partial x}  \log u$. Taking $x \rightarrow +\infty$, it follows from Lemma~\ref{lmmlimvar} that

\[ \lim_{x \rightarrow +\infty} \mathbb{E}_{\mu_{N,a}} [w(x)^2] = \mu^2 + \frac{\epsilon^2}{\kappa}\Gamma(0)\]

\noindent and that $\lim_{x \rightarrow +\infty} \mathbb{E}_{\mu_{N,a}}[w(x)]_x = \lim_{x \rightarrow +\infty} \mathbb{E}_{\mu_{N,a}}[u(x)] = 0$. From this, it follows that: 
  
\begin{equation}\label{eqmugam} 0 = - \frac{\kappa}{2}\left(\mu^2 + \frac{\epsilon^2}{\kappa}\Gamma(0)\right)+  \gamma \mu  - 1 + \frac{\epsilon^2}{2}\Gamma(0) \Rightarrow \mu^2 - \frac{2 \gamma}{\kappa} \mu + \frac{2}{\kappa} = 0.\end{equation}

\paragraph{The case $N > \sqrt{\frac{2}{\kappa}}$} From this, it follows directly that $\gamma  \geq \sqrt{2\kappa}$ (otherwise the equation has no solution for any $\mu \in \mathbb{R}$). Equation~\eqref{eqbdaprior} gives that if $N > \sqrt{\frac{2}{\kappa}}$, then   $\gamma  \leq \sqrt{2\kappa}$, so that:  

\[  N > \sqrt{\frac{2}{\kappa}} \Rightarrow \gamma  = \sqrt{2\kappa} \Rightarrow \mu = \sqrt{\frac{2}{\kappa}}.\]

\paragraph{The case $N < \sqrt{\frac{2}{\kappa}}$} Firstly, from Lemma~\ref{lmmgotlim}, together with the fact that $\gamma =  \lim_{t  \rightarrow +\infty}\frac{g(t)}{t}$, gives that

\[ \gamma \leq \frac{\kappa N}{2} + \frac{1}{N}.\]

\noindent Secondly, $\mu = \lim_{x \rightarrow +\infty} \left (\lim_{t \rightarrow +\infty} \mathbb{E}[w(t,x)]\right) \leq N$ and, from~\eqref{eqmugam},

 \[ \gamma = \frac{\kappa \mu}{2} + \frac{1}{\mu}.\]
 
 \noindent Since $f(\mu) = \frac{\kappa\mu}{2} + \frac{1}{\mu}$ is a  decreasing function in  $\mu$ for $\mu \in \left [0,  \sqrt{\frac{2}{\kappa}} \right )$, with $f(0) = +\infty$, it follows that  in the range $N < \sqrt{\frac{2}{\kappa}}$, $\mu = N$. The result follows. 
\qed  

\subsection{Correspondence between Correlated Noise and Standard Wiener Noise} Let $u$ solve~\eqref{eqspace} and let $u^*$ solve the same equation with initial condition $u^*(0,.) \equiv 1$. Let $\widehat{u} = \frac{u}{u^*}$. Then an application of Itô's formula gives:

\[ \left\{ \begin{array}{l} \widehat{u}_t = \frac{\kappa}{2}\widehat{u}_{xx} + \kappa (\log u^*)_x \widehat{u}_x + u^* \widehat{u}(1 - \widehat{u}) \\ \widehat{u}(0,x) = {\bf 1}_{(-\infty,-\frac{1}{N}\log 2]} + e^{-Nx}{\bf 1}_{(-\frac{1}{N}\log 2, +\infty)} \end{array}\right. \]

\noindent The foregoing shows that for $\epsilon > 0$, $\mathbb{E}[u^*(t,x)] < 1$ for $t > 0$, which should slow down the speed of the travelling front. The slow-down corresponding to the fact that $\overline{u^*} < 1$ is compensated by the additional drift term.

\section{Conclusion} The way that the perturbation $\epsilon u dW$ affects the speed of the travelling front for the KPP equation is already known from the work of Øksendal, Våge and Zhao~\cite{OkVaZh}, although this article presents a different proof which presents sharper results and is hopefully more transparent.

The main aim of the paper is to compare the effect of the perturbation $\epsilon u dW$ with $\epsilon u d\zeta$, where $W$ is a standard Wiener process (no spatial dependence), while the spatial covariance of $\zeta$ is integrable. 

The main result is that the spatial variation in the noise introduces a drift which compensates the slow-down that results from introducing noise, so that in this situation, the speed of the travelling front remains the same.

This observation raises a number of interesting questions, when these results are compared with those of Mueller et. al., where noise of the form $\epsilon u^{1/2} \xi$ results in a markedly {\em greater} slow-down of the travelling front. A coherent and unified picture of how noise affects the speed of the travelling front is an open problem.

 \setlength{\baselineskip}{2ex}

\end{document}